\documentclass{amsart}

\usepackage{hyperref, amssymb, amsmath, amscd}
\usepackage[all]{xy}
\usepackage{paralist, enumerate}

\allowdisplaybreaks[4]

\theoremstyle{plain}
\newtheorem{theorem}{Theorem}[section]
\newtheorem{proposition}[theorem]{Proposition}
\newtheorem{lemma}[theorem]{Lemma}
\newtheorem{corollary}[theorem]{Corollary}

\theoremstyle{definition}
\newtheorem{definition}[theorem]{Definition}
\newtheorem{remark}[theorem]{Remark}
\newtheorem{notation}[theorem]{Notation}
\newtheorem{example}[theorem]{Example}

\newcommand{\CC}{\mathbb{C}}
\newcommand{\KK}{\mathbb{K}}
\newcommand{\MM}{\mathbb{M}}
\newcommand{\PP}{\mathbb{P}}
\newcommand{\QQ}{\mathbb{Q}}
\newcommand{\SSS}{\mathbb{S}}
\newcommand{\VV}{\mathbb{V}}
\newcommand{\WW}{\mathbb{W}}
\newcommand{\ZZ}{\mathbb{Z}}

\newcommand{\uu}{\mathbf{u}}
\newcommand{\vv}{\mathbf{v}}
\newcommand{\ww}{\mathbf{w}}
\newcommand{\xx}{\mathbf{x}}
\newcommand{\yy}{\mathbf{y}}

\newcommand{\F}{\mathcal{F}}
\newcommand{\Hsh}{\mathcal{H}}
\newcommand{\I}{\mathcal{I}}
\newcommand{\Osh}{\mathcal{O}}
\newcommand{\Psh}{\mathcal{P}}
\newcommand{\X}{\mathcal{X}}

\newcommand{\A}{\mathfrak{A}}
\newcommand{\m}{\mathfrak{m}}
\newcommand{\p}{\mathfrak{p}}
\newcommand{\rid}{\mathfrak{r}}

\newcommand{\pr}{\mathrm{pr}}

\newcommand{\trans}{{\small \mathrm{T}}}

\newcommand{\der}{\operatorname{\Delta}}
\newcommand{\coh}{\operatorname{H}}
\newcommand{\depth}{\operatorname{depth}}
\newcommand{\GL}{\operatorname{GL}}
\newcommand{\hilb}{\operatorname{\mathcal{H}}}
\newcommand{\height}{\operatorname{ht}}

\newcommand{\lcm}{\operatorname{lcm}}
\newcommand{\lead}{\operatorname{lead}}
\newcommand{\length}{\operatorname{length}}
\newcommand{\rad}{\operatorname{rad}}
\newcommand{\red}{\operatorname{red}}
\newcommand{\reg}{\operatorname{reg}}
\newcommand{\rem}{\operatorname{rem}}
\newcommand{\rin}{\operatorname{rin}}
\newcommand{\room}{\operatorname{span}}

\newcommand{\lra}{\longrightarrow}
\newcommand{\tensor}{\otimes}

\begin{document}

\setdefaultleftmargin{21pt}{}{}{}{}{}

\title[On the Hilbert Function of a Finite Scheme]
{The Relevant Domain of the Hilbert Function of a Finite Multiprojective Scheme}

\author{Mario Maican}
\address{Institute of Mathematics of the Romanian Academy, Calea Grivitei 21, Bucharest 010702, Romania}

\email{maican@imar.ro}

\begin{abstract}
Let $X$ be a zero-dimensional scheme contained in a multiprojective space.
Let $s_i$ be the length of the projection of $X$ onto the $i$-th component
of the multiprojective space.
A result of Van~Tuyl states that
the Hilbert function of $X$, in the case when $X$ is reduced,
is completely determined by its restriction to the product of the intervals $[0, s_i - 1]$.
We prove that the same is also true for non-reduced schemes $X$.
\end{abstract}

\subjclass[2010]{Primary 13D40, 14N05, 14F06}
\keywords{Hilbert function, Zero-dimensional scheme, Multiprojective scheme}

\maketitle


\section{Introduction}
\label{introduction}

\noindent
Let $\VV = \PP^{n_1} \times \dots \times \PP^{n_q}$ be a multiprojective space over a field $\KK$.
Here $q$ and $n_1^{}, \dots, n_q^{}$ are positive integers.
The coordinate ring of $\VV$ is the $\ZZ^q$-graded algebra
\[
\SSS = \KK[x_{ij}^{} \mid 1 \le i \le q, \ 0 \le j \le n_i^{}].
\]
We have $\deg (x_{ij}^{}) = e_i^{}$,
where $e_i^{} \in \ZZ^q$ is the $i$-th basis element.
Let $M$ be a finitely generated $\ZZ^q$-graded $\SSS$-module.
The \emph{Hilbert function} $\hilb_M^{} \colon \ZZ^q \to \ZZ$ of $M$ is defined by
$\hilb_M^{}(a) =  \dim_\KK^{}(M_a^{})$.
Let $X \subset \VV$ be a zero-dimensional subscheme
and let $I(X) \subset \SSS$ be the ideal generated by the $\ZZ^q$-homogeneous forms in $\SSS$ that vanish on $X$.
The Hilbert function $\hilb_X^{}$ of $X$ is defined to be the Hilbert function of $\SSS / I(X)$.

The exploration of the Hilbert functions in the multiprojective setting is a natural extension
of the rich theory of Hilbert functions of zero-dimensional subschemes of $\PP^n$.
The simplest case, when $\VV = \PP^1 \times \PP^1$, was first investigated
by Giuffrida et al.\ in \cite{giuffrida_postulation}.
This exploration was then continued by many authors.
The case of $\PP^1 \times \PP^1$ remains the most assiduously studied case, see
\cite{carlini_fat}, \cite{projective_unexpected}, \cite{favacchio_three}, \cite{giuffrida_scheme},
\cite{guardo_quadric}, \cite{guardo_hilbert}, \cite{guardo_minimal}, \cite{guardo_points}, \cite{guardo_classifying},
\cite{marino_conductor}, \cite{marino_quadric}, \cite{tuyl_hilbert}.
For other ambient spaces $\VV$ we refer to
\cite{projective_unexpected}, \cite{favacchio_points}, \cite{favacchio_multiprojective},
\cite{tuyl_border}, \cite{tuyl_hilbert}.
The theory now follows three broad directions of development:\
the Hilbert functions of sets of points, as in
\cite{favacchio_points}, \cite{guardo_points}, \cite{guardo_classifying},
\cite{marino_conductor}, \cite{marino_quadric};
the Hilbert functions of sets of fat points, as in
\cite{carlini_fat}, \cite{projective_unexpected}, \cite{favacchio_three},
\cite{guardo_quadric}, \cite{guardo_hilbert}, \cite{guardo_minimal}, \cite{sidman_regularity};
and the Hilbert functions of ACM schemes, as in
\cite{favacchio_points}, \cite{favacchio_multiprojective},
\cite{guardo_quadric}, \cite{guardo_hilbert}, \cite{guardo_points}, \cite{guardo_classifying},
\cite{marino_conductor}, \cite{marino_quadric}, \cite{tuyl_hilbert}.

In connection with the first and the third directions of development,
we mention two fundamental results belonging to Van~Tuyl.
According to \cite{tuyl_border}, if $X$ is reduced and $\KK$ is algebraically closed,
then $\hilb_X^{}$ is uniquely determined by its restriction to
a rectangular region of the form $R = [0, r_1^{}] \times \cdots \times [0, r_q^{}] \subset \ZZ^q$.
More precisely, $R$ is a relevant domain for $\hilb_X^{}$
in the sense of Definition~\ref{relevance} and Lemma~\ref{relevant}.
Our first achievement in this paper is the generalization of this result
to the case of an arbitrary zero-dimensional subscheme $X \subset \VV$
and an arbitrary ground field $\KK$.
See Theorem~\ref{main}.
According to \cite{tuyl_hilbert}, if $\KK$ is algebraically closed,
then the functions $\hilb_X^{}$,
where $X$ runs through the zero-dimensional reduced ACM subschemes of $\VV$,
are precisely the functions $\Hsh$ whose difference $\der \Hsh$ (see formula \eqref{derivation})
is the Hilbert function of an artinian $\ZZ^q$-graded quotient of $\SSS / (x_{1 0}^{}, \dots, x_{q 0}^{})$.
Our second achievement in this paper is the generalization of this result
to the case of an arbitrary zero-dimensional subscheme $X \subset \VV$
and an arbitrary infinite ground field $\KK$.
See Theorem~\ref{artinian}.

We give three proofs to Theorem~\ref{main}.
The first proof consists of comparing the cohomology of the twists of $\I_X^{}$
with the cohomology of the twists of the ideal sheaf of $X$ in a smaller ambient space $W_i^{} \subset \VV$.
Here $W_i^{}$ is obtained from $\VV$ by replacing $\PP^{n_i}$ with the projection of $X$ onto $\PP^{n_i}$.
See Lemma~\ref{reduction}.
The second proof, located in section~\ref{constraints}, applies only in the case when
$\KK$ is algebraically closed, $\VV = (\PP^1)^q$ and $X$ is ACM or sub-ACM
(meaning $\depth(\SSS / I(X)) = q - 1$).
The technique we use draws on the technique of Giuffrida et al.,
who dealt with the case when $\VV = \PP^1 \times \PP^1$.
The key ingredient here are the constraints
satisfied by the Hilbert function of an ACM or sub-ACM scheme $X$
(Propositions~\ref{constraint_1} and \ref{constraint_2}).
The third proof of Theorem~\ref{main} is located in section~\ref{end}
and applies only in the case when $\KK$ is infinite and $\VV = (\PP^1)^q$.
It is based on Macaulay's theorem (Theorem~\ref{macaulay})
and on a vanishing criterion for the difference $\der \hilb_{\SSS / J}$
of the Hilbert function of the quotient of $\SSS$ by a monomial ideal (Proposition~\ref{vanishing}).
We think that approaching Theorem~\ref{main} from three different angles
provides a clearer picture of the subtleties that arise in the study of multiprojective
Hilbert functions.

An important consequence of Theorem~\ref{main} is
an upper estimate on the regularity index of $\SSS / I(X)$, regarded as a $\ZZ$-graded $\SSS$-module,
in terms of the regularity indices of the projections of $X$ onto the components $\PP^{n_i}$ of $\VV$.
See Corollary~\ref{estimate}.

Van~Tuyl's method for proving his version of Theorem~\ref{artinian} consists
of finding a regular sequence $\{ u_1^{}, \dots, u_q^{} \}$ for $\SSS / I(X)$, as in Proposition~\ref{regular_1}.
We adapt Van~Tuyl's argument
to the case when $X \subset \VV$ is an arbitrary zero-dimensional subscheme
and $\KK$ is an arbitrary infinite field.

In this paper we also consider quasi-rectangular domains,
i.e.\ finite unions of rectangular domains, that are relevant to $\hilb_X^{}$
in the sense of Definition~\ref{relevance}.
The third achievement of this paper is Proposition~\ref{procedure},
which gives sufficient conditions for the existence of quasi-rectangular relevant domanins
that are strictly contained in $R$.
The problem of describing all quasi-rectangular domains $Q \subset R$ that are relevant
for $\hilb_X^{}$ remains open.
An important class of schemes $X$ for which this problem has been settled
is the class of ACM subschemes of $(\PP^1)^q$.
See Corollary~\ref{corners}.

We now present the outline of the paper.
In section~\ref{domains} we gather a few elementary facts about relevant domains.
In section~\ref{generalities} we collect a few well-known facts about
Hilbert functions of finite subschemes of $\PP^n$.
These facts will be needed in the proof of our first main theorem,
concerning the rectangular relevant domain,
to whom section~\ref{rectangular} is devoted.
In section~\ref{regular}, whose role is to prepare the ground for the next two sections,
we construct regular sequences for $\SSS / I(X)$ and for $I(X)$
in the case when $X$ is ACM or sub-ACM.
Section~\ref{ACM} contains our second main theorem, concerning ACM schemes.
In section~\ref{constraints} we combine the results of section~\ref{regular}
with Lemma~\ref{koszul} in order to obtain inequalities involving the partial difference functions
of $\hilb_X^{}$ and $\hilb_{I(X)}^{}$.
As an application, we obtain our second proof of Theorem~\ref{main}.
In section~\ref{formula} we find a formula for $\der \hilb_{\SSS / J}^{}$,
where $J$ is a monomial ideal.
This leads us to our vanishing criterion for $\der \hilb_{\SSS / J}^{}$.
Section~\ref{end} contains our third proof of Theorem~\ref{main}
and our procedure for detecting quasi-rectangular relevant domains.


\section{Relevant domains}
\label{domains}

\noindent
Let $q \ge 1$ be an integer.
We introduce a partial order on $\ZZ^q$ as follows:\
given $a = (a_1^{}, \dots, a_q^{})$ and $b = (b_1^{}, \dots, b_q^{})$ in $\ZZ^q$,
we say that $a \le b$ if $a_i^{} \le b_i^{}$ for all indices $i \in \{ 1, \dots, q \}$.
Let $e_i^{} = (0, \dots, 1, \dots, 0)$ be the element of $\ZZ^q$
that has entry $1$ on position $i$ and entries $0$ elsewhere.
Let $\F \colon \ZZ^q \to \ZZ$ be a function.
We introduce the \emph{difference} function $\der \F \colon \ZZ^q \to \ZZ$ by the formula
\begin{equation}
\label{derivation}
\der \F(a) = \F(a) +
\sum_{1 \le p \le q} (-1)^p \sum_{1 \le i_1 < \dots < i_p \le q} \F(a - e_{i_1}^{} - \dots - e_{i_p}^{}).
\end{equation}
In this paper we shall only consider functions $\F$
that vanish on the complement of the positive quadrant $\ZZ_+^q = \{ a \in \ZZ^q \mid a \ge 0 \}$
because we are chiefly interested in Hilbert functions of $\ZZ^q$-homogeneous ideals in $\SSS$.
For such functions we can recover $\F$ from $\der \F$ by means of the formula
\begin{equation}
\label{integration}
\F(a) = \sum_{0 \le b \le a} \der \F(b) \quad \text{for all} \quad a \in \ZZ_+^q.
\end{equation}
Given $r, s \in \ZZ^q$ such that $r \le s$, we write $[r, s] = \{ a \in \ZZ^q \mid r \le a \le s \}$.
A \emph{rectangular domain} in $\ZZ^q$ has the form $R = [0, r]$, for some $r \in \ZZ_+^q$.
A \emph{quasi-rectangular domain} $Q \subset \ZZ^q$ is a finite union of rectangular domains.
The \emph{boundary} $B_Q^{}$ of $Q$ is defined to be the boundary of $Q$ inside $\ZZ_+^q$:
\[
B_Q^{} = \{ a \in Q \mid a + e_{i_1}^{} + \dots + e_{i_p}^{} \notin Q \ \text{for some indices} \ 1 \le i_1 < \dots < i_p \le q \}.
\]
In particular, for $R = [0, r]$,
$
B_R^{} = \{ a \in R \mid a_i^{} = r_i^{} \ \text{for some index} \ 1 \le i \le q \}.
$

\begin{definition}
\label{relevance}
Under the above notations, a quasi-rectangular domain $Q$ is said to be \emph{relevant} to $\F$
if $\der \F(a) = 0$ for all $a \in \ZZ^q \setminus Q$.
\end{definition}

\begin{lemma}
\label{relevant}
A rectangular domain $[0, r]$ is relevant to $\F$ if and only if
for every index $i \in \{1, \dots, q \}$ and for every $a \in \ZZ^q$ such that $a_i^{} \ge r_i^{}$
we have the equation $\F(a) = \F(a_1^{}, \dots, r_i^{}, \dots, a_q^{})$.
\end{lemma}

\begin{proof}
Assume that $R = [0, r]$ is relevant to $\F$ and choose $a \in \ZZ_+^q$ such that $a_i^{} \ge r_i^{}$.
Formula~\eqref{integration} can be rewritten as
\[
\F(a) = \sum_{0 \le b_1 \le a_1} \!\!\! \cdots \!\!\! \sum_{0 \le b_i \le r_i} \!\!\! \cdots \!\!\! \sum_{0 \le b_q \le a_q} \!\!\! \der \F(b)
+ \sum_{0 \le b_1 \le a_1} \!\!\! \cdots \!\!\! \sum_{r_i < b_i \le a_i} \!\!\! \cdots \!\!\! \sum_{0 \le b_q \le a_q} \!\!\! \der \F(b).
\]
The first summation equals $\F(a_1^{}, \dots, r_i^{}, \dots, a_q^{})$, again by virtue of formula~\eqref{integration}.
The second summation vanishes because $\der \F(b) = 0$ if $b$ lies outside $R$.

Conversely, assume that for every index $i \in \{1, \dots, q \}$
and for every $a \in \ZZ^q$ such that $a_i^{} \ge r_i^{}$
we have the equation $\F(a) = \F(a_1^{}, \dots, r_i^{}, \dots, a_q^{})$.
This is equivalent to saying that for every index $i \in \{1, \dots, q \}$
and for every $a \in \ZZ^q$ such that $a_i^{} > r_i^{}$
we have the equation $\F(a) = \F(a - e_i^{})$.
Choose $a \in \ZZ_+^q \setminus R$.
There is an index $i$ such that $a_i^{} > r_i^{}$.
Formula~\eqref{derivation} can be rewritten as
\begin{multline*}
\der \F(a) = \big(\F(a) - \F(a - e_i^{}) \big) \\
+ \sum_{1 \le p \le q - 1} (-1)^p \!\!\! \sum_{\substack{1 \le j_1 < \dots < j_p \le q \\ j_1, \dots, j_p \neq i}}
\big(\F(a - e_{j_1}^{} - \dots - e_{j_p}^{}) - \F(a - e_i^{} - e_{j_1}^{} - \dots - e_{j_p}^{}) \big).
\end{multline*}
All terms in parentheses vanish, hence $\der \F(a) = 0$.
Thus, $R$ is relevant to $\F$.
\end{proof}

\begin{lemma}
\label{constant}
Assume that the rectangular domain $[0, r]$ is relevant to $\F$.
We claim that $\F(a) = \F(r)$ for all $a \in \ZZ^q$ such that $a \ge r$.
\end{lemma}

\begin{proof}
Applying Lemma~\ref{relevant} repeatedly, we obtain the equations
\[
\F(a) = \F(r_1^{}, a_2^{}, \dots, a_q^{}) = \F(r_1{}, r_2^{}, a_3^{}, \dots, a_q^{}) = \dots = \F(r). \qedhere
\]
\end{proof}

\begin{remark}
\label{boundary}
The above lemmas show that, if $R$ is relevant to $\F$,
then $\F|_{B_R}^{}$ determines the function $\F$ on the complement of $R$.
The same is true for a relevant quasi-rectangular domain $Q$.
Take for instance $Q = [0, s] \setminus [r, s]$ in $\ZZ^2$, where $0 < r_1^{} < s_1^{}$, $0 < r_2^{} < s_2^{}$,
and take $a \in [r, s]$.
We have the equation
\[
\F(a) = \F(a_1^{}, r_2^{} - 1) + \F(r_1^{} - 1, a_2^{}) - \F(r_1^{} - 1, r_2^{} - 1)
\]
and $(a_1^{}, r_2^{} - 1)$, $(r_1^{} - 1, a_2^{})$, respectively, $(r_1^{} - 1, r_2^{} - 1)$ lie on $B_Q^{}$.
\end{remark}

\noindent
Given $a \in \ZZ_+^q$ we write $| a | = a_1^{} + \dots + a_q^{}$.
Consider a function $\F \colon \ZZ^q \to \ZZ$ that vanishes on the complement of $\ZZ_+^q$.
We construct a new function $\overline{\F} \colon \ZZ_+ \to \ZZ$ by the formula
\[
\overline{\F}(d) = \sum_{\substack{a \in \ZZ_+^q \\ | a | = d}} \F(a).
\]

\begin{lemma}
\label{polynomial}
Let $\F \colon \ZZ^q \to \ZZ$ be a function that vanishes on the complement of $\ZZ_+^q$.
Assume that the rectangular domain $[0, r]$ is relevant to $\F$.
We claim that the restriction of $\overline{\F}$ to $[| r |, \infty)$ is a polynomial function in $d$
with rational coefficients and with dominant term $\F(r) d^{q - 1} / (q - 1)!$
\end{lemma}

\begin{proof}
Write $R = [0, r]$ and $r = (r_1^{}, \dots, r_q^{})$.
Given $b \in B_R^{}$ there is an integer $p = p_b^{} \in \{ 1, \dots, q \}$
and there are indices $1 \le i_1^{} < \dots < i_p^{} \le q$ such that
$b_{i_1}^{} = r_{i_1}^{}$, $\dots$, $b_{i_p}^{} = r_{i_p}^{}$
and $b_i^{} < r_i^{}$ for $i \in \{ 1, \dots, q \} \setminus \{ i_1^{}, \dots, i_p^{} \}$.
We consider the set
\[
A(b) = \{ a \in \ZZ_+^q \mid a_{i_1}^{} \ge r_{i_1}^{}, \dots, a_{i_p}^{} \ge r_{i_p}^{}, \,
a_i^{} = b_i^{} \ \text{for} \ i \in \{ 1, \dots, q \} \setminus \{ i_1^{}, \dots, i_p^{} \} \}.
\]
According to Lemma~\ref{relevant}, $\F(a) = \F(b)$ for all $a \in A(b)$.
For $d \ge | b |$ we also consider the set
\[
A_d^{}(b) = \{ a \in A(b) \mid | a | = d \}.
\quad \text{We notice that} \quad
| A_d^{}(b) | = {d - | b | + p - 1 \choose p - 1}.
\]
Assume now that $d \ge | r |$.
We have the decomposition
\[
\{ a \in \ZZ_+^q \mid | a | = d \} = \bigsqcup_{b \in B_R} A_d^{}(b).
\]
We now calculate:
\begin{align*}
\overline{\F}(d) & = \sum_{\substack{a \in \ZZ_+^q \\ | a | = d}} \F(a)
= \sum_{b \in B_R} \sum_{a \in A_d(b)} \F(a)
= \sum_{b \in B_R} \F(b) | A_d^{}(b) | \\
& = \sum_{b \in B_R} \F(b) {d - | b | + p_b^{} - 1 \choose p_b^{} - 1}.
\end{align*}
This is a polynomial expression in $d$ with rational coefficients.
The claim about the dominant term follows from the fact that $p_r^{} = q$
while $p_b^{} < q$ for $b \in B_R^{} \setminus \{ r \}$.
\end{proof}

\begin{definition}
\label{poincare}
In the context of the above lemma,
the polynomial $\Psh(d) \in \QQ[d]$ satisfying the condition $\Psh(d) = \overline{\F}(d)$ for $d \ge | r |$
will be called the \emph{Poincar\`e polynomial} associated to $\F$.
\end{definition}


\section{Generalities concerning Hilbert functions}
\label{generalities}

\noindent
Let $\VV = \PP^{n_1} \times \dots \times \PP^{n_q}$ be a multiprojective space over a field $\KK$.
Let $X \subset \VV$ be a zero-dimensional subscheme with ideal sheaf $\I_X^{}$ and structure sheaf $\Osh_X^{}$.
Denote $\length(X) = \dim_\KK^{} \coh^0(\Osh_X^{})$.
Choose $a \in \ZZ^q$.
From the short exact sequence
\[
0 \lra \I_X^{} \lra \Osh_\VV^{} \lra \Osh_X^{} \lra 0
\]
we obtain the exact sequence in cohomology
\[
0 \lra I(X)_a^{} \lra \SSS_a^{} \lra \coh^0(\Osh_X^{}) \lra \coh^1(\I_X^{}(a)) \lra \coh^1(\Osh_\VV^{}(a)).
\]
The group on the right vanishes if $a \ge 0$.
Thus,
\begin{equation}
\label{cohomology}
\hilb_X^{}(a) = \length(X) - \dim_\KK^{} \coh^1(\I_X^{}(a)) \quad \text{if} \quad a \ge 0.
\end{equation}

\begin{lemma}
\label{length}
Assume that $[0, r]$ is a relevant domain for $\hilb_X^{}$.
We claim that $\hilb_X^{}(a) = \length(X)$ for all $a \ge r$.
\end{lemma}

\begin{proof}
According to \cite[Theorem III.5.2]{hartshorne_algebraic},
$\coh^1(\I_X^{}(a))$ vanishes if $a_1^{}, \dots, a_q^{}$ are sufficiently large.
Consequently, in view of formula~\eqref{cohomology},
$\hilb_X^{}(a) = \length(X)$ if $a_1^{}, \dots, a_q^{}$ are sufficiently large.
We saw at Lemma~\ref{constant} that $\hilb_X^{}$ is constant on $[r, \infty)$.
We conclude that $\hilb_X^{}$ takes the value $\length(X)$ on $[r, \infty)$.
\end{proof}

\noindent
It is well-known that $\hilb_X^{}$ has a relevant domain when $\VV$ is a projective space.
The following proposition is a straightforward consequence of \cite[Proposition 1.1]{geramita_ideal}.

\begin{proposition}
\label{regularity_1}
Let $Z \subset \PP^n$ be a zero-dimensional subscheme.
We claim that there is an integer $r \ge 0$
such that $\hilb_Z^{}$ increases on the interval $[0, r]$ and is constant on the interval $[r, \infty)$.
Thus, $[0, r]$ is the smallest relevant domain for $\hilb_X^{}$.
\end{proposition}

\begin{notation}
\label{rin}
The integer $r = \rin(Z)$ is known as the \emph{regularity index} of $Z$.
\end{notation}

\noindent
In the following proposition we collect several well-known properties of the regularity index.
For the convenience of the reader we include their proofs.

\begin{proposition}
\label{regularity_2}
Let $Z \subset \PP^n$ be a zero-dimensional subscheme of length $s$ and regularity index $r$.
We make the following claims:
\begin{enumerate}
\item[\emph{(i)}]
$\hilb_Z^{}(a) = s$ for $a \ge r$;
\item[\emph{(ii)}]
$\coh^1(\I_Z^{}(a)) = \{ 0 \}$ for $a \ge r$;
\item[\emph{(iii)}]
$\coh^m(\I_Z^{}(a)) = \{ 0 \}$ for $m \ge 2$ and $a \ge -n$;
\item[\emph{(iv)}]
$r \le s - 1$;
\item[\emph{(v)}]
$r = s - 1$ if $n = 1$;
\item[\emph{(vi)}]
$r + 1 = \reg(Z)$, the Castelnuovo-Mumford regularity of $Z$.
\end{enumerate}
\end{proposition}

\begin{proof}
(i) We apply Proposition~\ref{regularity_1} and Lemma~\ref{length}.
(ii) In light of formula~\eqref{cohomology},
we have the equations $\dim_\KK^{} \coh^1(\I_Z^{}(a)) = s - \hilb_Z^{}(a) = 0$ for $a \ge r$.
(iii) From the exact sequence
\[
0 \lra \I_Z^{} \lra \Osh_{\PP^n}^{} \lra \Osh_Z^{} \lra 0
\]
we obtain the exact sequence
\[
\coh^{m - 1}(\Osh_Z^{}) \lra \coh^m(\I_Z^{}(a)) \lra \coh^m(\Osh_{\PP^n}^{}(a)).
\]
The group on the left vanishes because $\Osh_Z^{}$ has support of dimension zero.
The group on the right vanishes for $a \ge -n$.
Thus, the group in the middle also vanishes.
(iv) According to Proposition~\ref{regularity_1}, $\hilb_Z^{}$ increases on the interval $[0, r]$.
By definition, $\hilb_Z^{}(0) = 1$, hence $\hilb_Z^{}(r) \ge r + 1$, and hence $s \ge r + 1$.
(v) For $0 \le a \le s - 1$ we have $\hilb_Z^{}(a) = a + 1$
because there are no forms of degree $a$ vanishing on a subscheme $Z \subset \PP^1$ of length $s$.
Thus, $\hilb_Z^{}$ increases on the interval $[0, s - 1]$ forcing the inequality $s - 1 \le r$.
The reverse inequality was obtained at claim (iv) above.
(vi) From claims (ii) and (iii) we see that $\coh^m(\I_Z^{}(r + 1 - m)) = \{ 0 \}$ if $m \ge 1$.
From the definition of $r$ and from formula~\eqref{cohomology}
it follows that $\coh^1(\I_Z^{}(r - 1)) \neq \{ 0 \}$.
\end{proof}

\begin{notation}
\label{complex}
Let $S$ be a $\KK$-algebra and let $M$ be an $S$-module.
Consider elements $v_1^{}, \dots, v_p^{} \in S$.
We denote by $\{ \epsilon_1^{}, \dots, \epsilon_p^{} \}$
the standard basis of the $\KK$-vector space $E = \KK^p$.
Consider the element $\vv = \epsilon_1^{} \tensor v_1^{} + \dots + \epsilon_p^{} \tensor v_p^{} \in E \tensor_\KK^{} S$.
The sequence
\[
0 \to M \overset{\cdot \vv}{\lra} E \tensor M \to \dots
\wedge^k E \tensor M \overset{\cdot \vv}{\lra} \wedge^{k + 1} E \tensor M \to \dots
\wedge^{p - 1} E \tensor M \overset {\cdot \vv}{\lra} \wedge^p E \tensor M
\]
is the \emph{Koszul complex} associated to $v_1^{}, \dots, v_p^{}$ and $M$,
denoted $K(v_1^{}, \dots, v_p^{}) \tensor M$.
\end{notation}

\begin{lemma}
\label{koszul}
Let $M$ be a $\ZZ^q$-graded $\SSS$-module.
Assume that $\{ u_1^{}, \dots, u_q^{} \}$ is an $M$-regular sequence
with $u_i^{} \in \room \{ x_{ij}^{} \mid 0 \le j \le n_i^{} \}$.
We claim that $\der \hilb_M^{}$ is the Hilbert function of $M / (u_1^{}, \dots, u_q^{}) M$.
\end{lemma}

\begin{proof}
We denote by $\{ \epsilon_1^{}, \dots, \epsilon_q^{} \}$
the standard basis of the $\KK$-vector space $E = \KK^q$.
For each $k \in \{ 1, \dots, q \}$ we endow $\wedge^k E \tensor_\KK^{} M$ with a $\ZZ^q$-grading as follows:
if $h \in M$ is $\ZZ^q$-homogeneous and $1 \le i_1^{} < \dots < i_k^{} \le q$, then
\[
\deg(\epsilon_{i_1}^{} \wedge \dots \wedge \epsilon_{i_k}^{} \tensor h)
= \deg(h) + \sum_{\substack{1 \le i \le q \\ i \notin \{ i_1^{}, \dots, i_k^{} \}}} e_i^{}.
\]
The Koszul complex $K(u_1^{}, \dots, u_q^{}) \tensor M$ (see Notation~\ref{complex}),
becomes a complex of $\ZZ^q$-graded $\SSS$-modules.
According to \cite[Corollary 17.5]{eisenbud_commutative}, $K(u_1^{}, \dots, u_q^{}) \tensor M$ is exact,
by virtue of the fact that $\{ u_1^{}, \dots, u_q^{} \}$ is $M$-regular.
We have an isomorphism $M \simeq \wedge^q E \tensor M$ of $\ZZ^q$-graded $\SSS$-modules
given by $h \mapsto \epsilon_1^{} \wedge \dots \wedge \epsilon_q^{} \tensor h$.
Thus, the cokernel of the last map in the Koszul complex is isomorphic to $M / (u_1^{}, \dots, u_q^{}) M$.
The lemma follows from the additivity of the Hilbert function on short exact sequences.
\end{proof}


\section{The rectangular relevant domain}
\label{rectangular}

\noindent
Let $\VV = \PP^{n_1} \times \dots \times \PP^{n_q}$ be a multiprojective space over a field $\KK$.
Let $X \subset \VV$ be a zero-dimensional subscheme.
Let $\pr_i^{} \colon \VV \to \PP^{n_i}$ be the projection onto the $i$-th component.
Let $Z_i^{} = \pr_i^{}(X)$ be the zero-dimensional subscheme of $\PP^{n_i}$ defined by the ideal
$I(Z_i^{}) = I(X) \cap \KK[x_{ij}^{} \mid 0 \le j \le n_i^{}]$.
Write $s_i^{} = \length(Z_i^{})$ and $r_i^{} = \rin(Z_i^{})$.
Let $W_i^{} = \pr_i^{-1}(Z_i^{})$ be the pull-back scheme,
i.e.\ the subscheme of $\VV$ defined by the ideal of $\SSS$ generated by $I(Z_i^{})$.
Noting that $X$ is a subscheme of $W_i^{}$,
we denote by $\I_{X, W_i}^{}$ the ideal sheaf of $X$ in $\Osh_{W_i}^{}$.

In this section we will focus on proving
that the domain $[0, r_1^{}] \times \dots \times [0, r_q^{}]$ is relevant to $\hilb_X^{}$.
As mentioned in the introduction, a similar result was proved by Van~Tuyl,
see Corollary~\ref{superfluous} below.
Van~Tuyl's approach was based on his version of Proposition~\ref{decomposition} from below.
Our approach is to replace the ambient space $\VV$ with $W_i^{}$.
Our main technical ingredient is the following lemma.

\begin{lemma}
\label{reduction}
Let $X \subset \VV$ be a zero-dimensional subscheme.
Fix an arbitrary index $i \in \{1, \dots, q \}$ and let $W_i^{}$ and $r_i^{}$ be as defined above.
Let $a \in \ZZ^q$ satisfy the conditions $a \ge 0$ and $a_i^{} \ge r_i^{}$.
We claim that
\[
\coh^1(\I_X^{}(a)) \simeq \coh^1(\I_{X, W_i}^{}(a - a_i^{} e_i^{})).
\]
\end{lemma}

\begin{proof}
By symmetry, we may assume that $i = 1$.
By virtue of the K\"unneth formula,
\[
\coh^m(\I_{W_1}^{}(a)) \simeq \bigoplus_{m_1 + \dots + m_q = m}
\coh^{m_1}(\I_{Z_1}^{}(a_1^{})) \tensor \coh^{m_2}(\Osh_{\PP^{n_2}}^{}(a_2^{})) \tensor \cdots \tensor
\coh^{m_q}(\Osh_{\PP^{n_q}}^{}(a_q^{})).
\]
By hypothesis, $a_1^{} \ge r_1^{}$, hence, in view of Proposition~\ref{regularity_2},
$\coh^{m_1}(\I_{Z_1}^{}(a_1^{})) = \{ 0 \}$ for $m_1^{} \ge 1$.
By hypothesis, $a_2^{} \ge 0, \dots, a_q^{} \ge 0$, hence
the higher cohomology groups of $\Osh_{\PP^{n_2}}^{}(a_2^{}), \dots, \Osh_{\PP^{n_q}}^{}(a_q^{})$ also vanish.
We deduce that $\coh^m(\I_{W_1}^{}(a)) = \{ 0 \}$ for $m \ge 1$.
From the short exact sequence
\[
0 \lra \I_{W_1}^{} \lra \I_X^{} \lra \I_{X, W_1}^{} \lra 0
\]
of sheaves on $\VV$ we obtain the long exact sequence
\[
\{ 0 \} = \coh^1(\I_{W_1}^{}(a)) \lra \coh^1(\I_X^{}(a)) \lra \coh^1(\I_{X, W_1}^{}(a)) \lra \coh^2(\I_{W_1}^{}(a)) = \{ 0 \}.
\]
The middle arrow must be an isomorphism.
The line bundle $\Osh_{\PP^{n_1}}^{}(a_1^{})$ is trivial on $Z_1^{}$
because $Z_1^{}$ is supported on finitely many points.
It follows that $\I_{X, W_1}^{}(a) \simeq \I_{X, W_1}^{}(0, a_2^{}, \dots, a_q^{})$.
We obtain the desired isomorphism
\[
\coh^1(\I_X^{}(a)) \simeq \coh^1(\I_{X, W_1}^{}(0, a_2^{}, \dots, a_q^{})). \qedhere
\]
\end{proof}

\begin{proposition}
\label{decomposition}
Let $X \subset \VV$ be a zero-dimensional subscheme.
Assume that the ground field $\KK$ is algebraically closed.
Fix an arbitrary index $i \in \{1, \dots, q \}$
and assume that $Z_i^{}$ is reduced, say $Z_i^{} = \{ P_1^{}, \dots, P_m^{} \}$.
For each index $k \in \{ 1, \dots, m \}$, we consider the scheme
\[
\WW_k^{} = \pr_i^{-1}(P_k^{}) \simeq \PP^{n_1} \times \dots \times \widehat{\PP}{}^{n_i} \times \dots \times \PP^{n_q}
\]
and we put $Y_k^{} = X \cap \WW_k^{}$.
Let $\hilb_{Y_k}^{}$ be the Hilbert function of $Y_k^{}$ as a subscheme of $\WW_k^{}$,
the latter being regarded as a multiprojective space.
Let $a \in \ZZ^q$ satisfy the conditions $a \ge 0$ and $a_i^{} \ge r_i^{}$.
We claim that
\[
\hilb_X^{}(a) = \sum_{1 \le k \le m} \hilb_{Y_k}^{}(a_1^{}, \dots, \widehat{a}_i^{}, \dots, a_q^{}).
\]
\end{proposition}

\begin{proof}
From the decomposition $W_i^{} = \WW_1^{} \sqcup \dots \sqcup \WW_m^{}$ we obtain the decomposition
\[
\coh^1(\I_{X, W_i}^{}(a - a_i^{} e_i^{}))
\simeq \bigoplus_{1 \le k \le m} \coh^1(\I_{Y_k, \WW_k}^{}(a_1^{}, \dots, \widehat{a}_i^{}, \dots, a_q^{})).
\]
Applying formula~\eqref{cohomology} and Lemma~\ref{reduction}, we calculate:
\begin{align*}
\hilb_X^{}(a) & = \length(X) - \dim_\KK^{} \coh^1(\I_X^{}(a)) \\
& = \length(X) - \dim_\KK^{} \coh^1(\I_{X, W_i}^{} (a - a_i^{} e_i^{})) \\
& = \sum_{1 \le k \le m} \length(Y_k^{}) -
\sum_{1 \le k \le m} \dim_\KK^{} \coh^1(\I_{Y_k, \WW_k}^{}(a_1^{}, \dots, \widehat{a}_i^{}, \dots, a_q^{})) \\
& = \sum_{1 \le k \le m} \hilb_{Y_k}^{}(a_1^{}, \dots, \widehat{a}_i^{}, \dots, a_q^{}). \qedhere
\end{align*}
\end{proof}

\noindent
The above result, in the particular case when $X$ is reduced and $a_i^{} \ge s_i^{} - 1$,
was obtained by Van~Tuyl using different methods.
Consult \cite[Proposition 4.2]{tuyl_border}.

\begin{definition}
\label{rem}
Let $X \subset \VV$ be a zero-dimensional subscheme.
For each $i \in \{ 1, \dots, q \}$, let $Z_i^{}$ be the projection of $X$ onto $\PP^{n_i}$.
The tuple
\[
\rem(X) = (\rin(Z_1^{}), \dots, \rin(Z_q^{}))
\]
(see Notation~\ref{rin})
will be called the \emph{regularity multiindex} of $X$.
\end{definition}

\begin{theorem}
\label{main}
Let $X \subset \PP^{n_1} \times \dots \times \PP^{n_q}$ be a zero-dimensional subscheme.
We claim that $R(X) = [0, \, \rem(X)]$ is the smallest rectangular relevant domain for $\hilb_X^{}$.
\end{theorem}

\begin{proof}
Write $\rem(X) = (r_1^{}, \dots, r_q^{})$.
Consider $a \in \ZZ_+^q$ satisfying the condition $a_i^{} \ge r_i^{}$ for some index $i \in \{ 1, \dots, q \}$.
According to Lemma~\ref{reduction},
the expression $\dim_\KK^{} \coh^1(\I_X^{}(b))$ remains constant as $b_i^{}$ varies in the interval $[r_i^{}, \infty)$
and $b_j^{}$ are nonnegative fixed integers for all indices $j \in \{ 1, \dots, q \} \setminus \{ i \}$.
Thus,
\[
\dim_\KK^{} \coh^1(\I_X^{}(a)) = \dim_\KK^{} \coh^1(\I_X^{}(a_1^{}, \dots, r_i^{}, \dots, a_q^{})).
\]
Applying formula~\eqref{cohomology}, we calculate:
\begin{align*}
\hilb_X^{}(a) & = \length(X) - \dim_\KK^{} \coh^1(\I_X^{}(a)) \\
& = \length(X) - \dim_\KK^{} \coh^1(\I_X^{}(a_1^{}, \dots, r_i^{}, \dots, a_q^{})) \\
& = \hilb_X^{}(a_1^{}, \dots, r_i^{}, \dots, a_q^{}).
\end{align*}
From Lemma~\ref{relevant} we deduce that $R(X)$ is relevant to $\hilb_X^{}$.
We cannot shrink $R(X)$ to a smaller rectangular relevant domain because, as seen at Proposition~\ref{regularity_1},
the function $\hilb_X^{}(a_i^{} e_i^{}) = \hilb_{Z_i}^{}(a_i^{})$ increases on the interval $[0, r_i^{}]$.
\end{proof}

\noindent
The above result, in the particular case when $\VV = \PP^1 \times \PP^1$ and $\KK$ is algebraically closed,
was obtained by Giuffrida et al.
Consult \cite[Remark 2.8 and Theorem 2.11]{giuffrida_postulation}.
Using different methods, Guardo and Van~Tuyl proved the above result in the particular case when
$\VV = \PP^1 \times \PP^1$, $\KK$ is algebraically closed and $X$ is a union of fat points.
Consult \cite[Corollary 3.4]{guardo_hilbert}.

\begin{proposition}
\label{regularity_3}
Let $X \subset \PP^{n_1} \times \dots \times \PP^{n_q}$ be a zero-dimensional subscheme
of length $s$ and regularity multiindex $r$.
For each index $i \in \{ 1, \dots, q \}$, let $s_i^{}$ be the length of the projection of $X$ onto $\PP^{n_i}$.
We make the following claims:
\begin{enumerate}
\item[\emph{(i)}]
$\hilb_X^{}(a) = s$ for $a \ge r$;
\item[\emph{(ii)}]
$\coh^1(\I_X^{}(a)) = \{ 0 \}$ for $a \ge r$;
\item[\emph{(iii)}]
$\coh^m(\I_X^{}(a)) = \{ 0 \}$ for $m \ge 2$ and $a \ge 0$;
\item[\emph{(iv)}]
$r \le (s_1^{} - 1, \dots, s_q^{} - 1)$;
\item[\emph{(v)}]
$r = (s_1^{} - 1, \dots, s_q^{} - 1)$ if $\VV = (\PP^1)^q$.
In particular, adopting the notations of Theorem~\ref{main},
$R(X) = [0, s_1^{} - 1] \times \dots \times [0, s_q^{} - 1]$.
\end{enumerate}
\end{proposition}

\begin{proof}
Claim (i) follows from Theorem~\ref{main} and Lemma~\ref{length}.
To prove claims (ii) and (iii) we argue precisely as in the proof of Proposition~\ref{regularity_2}(ii), respectively, (iii).
Claims (iv) and (v) follow from their counterparts at Proposition~\ref{regularity_2}.
\end{proof}

\noindent
Claim (i) of the above proposition,
in the particular case when $X$ is a union of fat points and $\KK$ is algebraically closed,
was obtained by Sidman and Van~Tuyl.
Consult \cite[Proposition 4.4]{sidman_regularity}.

\begin{corollary}
\label{superfluous}
Let $X \subset \PP^{n_1} \times \dots \times \PP^{n_q}$ be a zero-dimensional subscheme.
For each index $i \in \{ 1, \dots, q \}$, let $s_i^{}$ be the length of the projection of $X$ onto $\PP^{n_i}$.
We claim that the rectangular domain $[0, s_1^{} - 1] \times \dots \times [0, s_q^{} - 1]$
is relevant to $\hilb_X^{}$.
\end{corollary}

\noindent
The corollary follows from Theorem~\ref{main} and Proposition~\ref{regularity_3}(iv).
The above result, in the particular case when $X$ is reduced and $\KK$ is algebraically closed,
was obtained by Van~Tuyl.
Consult \cite[Proposition 4.6(ii) and Corollary 4.7]{tuyl_border}.

Let $M$ be a finitely generated $\ZZ^q$-graded $\SSS$-module.
Consider the canonical $\ZZ$-grading on $\SSS$ given by the degree of a polynomial.
Given $a \in \ZZ^q$, write $| a | = a_1^{} + \dots + a_q^{}$.
Then $M$ is also a $\ZZ$-graded $\SSS$-module with $M_d^{} = \oplus_{| a | = d}^{} M_a^{}$.
The Hilbert function of this module is the function $\overline{\hilb}_M^{} \colon \ZZ \to \ZZ$ given by
\[
\overline{\hilb}_M^{} (d) = \sum_{| a | = d} \hilb_M^{}(a).
\]
There exists a polynomial $\Psh_M^{}$ in one variable, with rational coefficients,
called the \emph{Hilbert-Poincar\`e polynomial} of $M$,
such that $\overline{\hilb}_M^{} (d) = \Psh_M^{}(d)$ for $d$ sufficiently large.
See \cite[Theorem 1.11]{eisenbud_commutative}.
The \emph{regularity index} of $M$, written $\rin(M)$, is the smallest integer
with the property that $\overline{\hilb}_M^{} = \Psh_M^{}$ on $[\rin(M), \infty)$.
Given a subscheme $X \subset \VV$, we put $\overline{\hilb}_X^{} = \overline{\hilb}_{\SSS / I(X)}^{}$,
$\Psh_X^{} = \Psh_{\SSS / I(X)}^{}$ and $\rin(X) = \rin(\SSS / I(X))$.

\begin{corollary}
\label{estimate}
Let $X \subset \VV$ be a zero-dimensional subscheme.
We claim that $\rin(X) \le | \rem(X) |$ (see Definition~\ref{rem})
and that the polynomial $\Psh_X^{}(d)$ has dominant term
$\length(X) d^{q - 1} / (q - 1)!$
\end{corollary}

\begin{proof}
The Poincar\`e polynomial associated to $\hilb_X^{}$ (see Definition~\ref{poincare})
coincides with $\Psh_X^{}$.
Thus, $\Psh_X^{} = \overline{\hilb}_X^{}$ on $[\rem(X), \, \infty)$,
hence $\rin(X) \le | \rem(X) |$.
According to Lemma~\ref{polynomial} and Proposition~\ref{regularity_3}(i),
$\Psh_X^{}(d)$ has dominant term
\[
\frac{\hilb_X^{}(\rem(X)) d^{q - 1}}{(q - 1)!} = \frac{\length(X) d^{q - 1}}{(q - 1)!}. \qedhere
\]
\end{proof}


\section{Regular sequences in the case of ACM and sub-ACM schemes}
\label{regular}

\noindent
In this section we assume that the ground field $\KK$ is infinite.
We denote by $\m$ the maximal ideal of $\SSS$ generated by all the variables.
We recall that the \emph{depth} of a $\ZZ^q$-graded $\SSS$-module $M$
is the maximal length of an $M$-regular sequence contained in $\m$.
Let $X \subset \VV$ be a zero-dimensional subscheme.
Since $\KK$ is infinite, there exists a non-constant $\ZZ^q$-homogeneous form
that does not vanish at any point of $\red(X)$.
This form is a non-zerodivisor of $\SSS / I(X)$.
We have the inequalities
\[
1 \le \depth(\SSS / I(X)) \le \dim(\SSS / I(X)) = q.
\]
We say that $X$ is \emph{arithmetically Cohen-Macaulay (ACM)}
if we have equality on the right.
If $q = 1$, then $X$ is automatically ACM.
We say that $X$ is \emph{sub-ACM} if $q \ge 2$ and $\depth(\SSS / I(X)) = q - 1$.

\begin{notation}
\label{U}
Throughout this section we shall employ the following notations:
\[
U = \room \{ x_{ij}^{} \mid 1 \le i \le q, \ 0 \le j \le n_i^{} \},
\]
\[
U_i^{} = \room \{ x_{ij}^{} \mid 0 \le j \le n_i^{} \},
\]
\[
U_i^0 = \{ u_i^{} \in U_i^{} \mid u_i^{} \ \text{does not vanish at any point of} \ \red(X) \}.
\]
\end{notation}

\begin{remark}
\label{U_0}
The spaces $U_i^0$ are non-empty and each $u_i^{} \in U_i^0$ is a non-zerodivisor for $\SSS / I(X)$.
Indeed, for each closed point $P \in \VV$, $I(P) \cap U_i^{}$ is a proper vector subspace of $U_i^{}$.
A vector space over an infinite field cannot be a finite union of proper subspaces.
The ideals $I(P)$ with $P \in \red(X)$ are the associated primes of $X$,
hence, by construction, $u_i^{}$ is a non-zerodivisor of $\SSS / I(X)$.
\end{remark}

\begin{remark}
\label{U_bis}
If $\KK$ is algebraically closed and $\red(X) = \{ P_1^{}, \dots, P_m^{} \}$,
then, for each index $k \in \{ 1, \dots, m \}$, there are vector subspaces $U_{ki}^{} \subset U_i^{}$ of codimension one
such that $I(P_k^{}) = (U_{ki}^{} \mid 1 \le i \le q)$.
We have $U_i^0 = U_i^{} \setminus \bigcup_{1 \le k \le m} U_{ki}^{}$.
\end{remark}

\noindent
In the case when $\KK$ is algebraically closed and $X$ is reduced and ACM,
Van~Tuyl proved that we can choose a regular sequence $\{ u_1^{}, \dots, u_q^{} \}$ for $\SSS / I(X)$ with $u_i^{} \in U_i^{}$.
Consult \cite[Proposition 3.2]{tuyl_hilbert}.
The aim of this section is to generalize this result to the case
when $X$ is an arbitrary zero-dimensional ACM subscheme
and $\KK$ is an arbitrary infinite field (see Proposition~\ref{regular_1}),
to obtain a version for sub-ACM schemes (see Proposition~\ref{regular_2}),
and to show that the $u_i^{}$ above can be chosen generically (see Proposition~\ref{regular_3}).
These results and their corollaries will be used in sections~\ref{ACM} and \ref{constraints}.

\begin{lemma}
\label{noether}
Assume that $\KK$ is infinite.
Let $\p \subset \KK[\xx, \yy]$ be a prime ideal,
where $\xx = (x_1^{}, \dots, x_m^{})$ and $\yy = (y_1^{}, \dots, y_n^{})$.
We claim that $\height(\p) \ge \height(\p \cap \KK[\xx]) + \height(\p \cap \KK[\yy])$.
\end{lemma}

\begin{proof}
Write $k = \dim \KK[\xx] / \p \cap \KK[\xx]$ and $l = \dim \KK[\yy] / \p \cap \KK[\yy]$.
Applying \cite[Corollary 13.4]{eisenbud_commutative},
we see that the claim is equivalent to the inequality $\dim \KK[\xx, \yy] / \p \le k + l$.
Applying the Noether normalization theorem \cite[Theorem 13.3]{eisenbud_commutative}
we deduce that there are linearly independent one-forms $u_1^{}, \dots, u_k^{} \in \KK[\xx]$,
respectively, $v_1^{}, \dots, v_l^{} \in \KK[\yy]$
such that $\KK[\xx] / \p \cap \KK[\xx]$ is integral over $\KK[\uu]$
and $\KK[\yy] / \p \cap \KK[\yy]$ is integral over $\KK[\vv]$.
Here $\uu = (u_1^{}, \dots, u_k^{})$ and $\vv = (v_1^{}, \dots, v_l^{})$.
Clearly, the extension of algebras $\KK[\uu, \vv] / \p \cap \KK[\uu, \vv] \subset \KK[\xx, \yy] / \p$ is integral,
hence the two algebras have the same dimension
(equal to the transcendence degree over $\KK$ of their fields of fractions,
see \cite[Theorem A, p.\ 286]{eisenbud_commutative}).
Thus, $\dim \KK[\xx, \yy] / \p \le \dim \KK[\uu, \vv] = k + l$.
\end{proof}

\begin{lemma}
\label{radical}
Let $X \subset \VV$ be a zero-dimensional subscheme.
Choose $u_i^{} \in U_i^0$.
We claim that $(U_i^{}) \subset \rad((u_i^{}) + I(X))$.
\end{lemma}

\begin{proof}
Recall that $I(X) \cap \KK[x_{ij}^{} \mid 0 \le j \le n_i^{}]$
defines a zero-dimensional subscheme $Z_i^{} \subset \PP^{n_i}$.
By hypothesis, $u_i^{}$ does not vanish at any point of $\red(Z_i^{})$,
hence $(u_i^{}) + I(Z_i^{})$ defines the empty subscheme in $\PP^{n_i}$,
and hence
\[
(U_i^{}) \subset \rad((u_i^{}) + I(Z_i^{})) \subset \rad((u_i^{}) + I(X)). \qedhere
\]
\end{proof}

\begin{proposition}
\label{regular_1}
Assume that $\KK$ is infinite.
Let $X \subset \VV$ be a zero-dimensional ACM subscheme.
We claim that there are  $u_i^{} \in U_i^{}$
such that $\{ u_1^{}, \dots, u_q^{} \}$ is a regular sequence for $\SSS / I(X)$.
\end{proposition}

\begin{proof}
Write $\xx_i^{} = (x_{i0}^{}, \dots, x_{i n_i})$.
Performing induction on $k \in \{ 1, \dots, q \}$, we will construct a regular sequence
$\{ u_1^{}, \dots, u_k^{} \}$ for $\SSS / I(X)$ with $u_i^{} \in U_i^0$.
To start the induction, choose $u_1^{} \in U_1^0$ and recall Remark~\ref{U_0}.
For the induction step, assume that $k \in \{ 1, \dots, q - 1 \}$
and that $\{ u_1^{}, \dots, u_k^{} \}$ has already been constructed.
Write $J = (u_1^{}) + \dots + (u_k^{}) + I(X)$.
Choose a prime ideal $\p$ that is associated to $J$.
By hypothesis, $\SSS / I(X)$ is Cohen-Macaulay of dimension $q$,
hence $\SSS / J$ is Cohen-Macaulay of dimension $q - k$.
According to \cite[Corollary 18.14]{eisenbud_commutative}, $J$ is unmixed,
hence $\dim(\p) = q - k$,
and hence, $\height(\p) = n_1^{} + \dots + n_q^{} + k$.
According to Lemma~\ref{radical}, $(U_1^{}, \dots, U_k^{})$ lies in $\rad(J)$,
so it is contained in $\p$.
We claim that $\height(\p \cap \KK[\xx_i^{}]) = n_i^{}$ for each index $i \in \{ k + 1, \dots, q \}$.
Indeed, the inequalities $\height(\p \cap \KK[\xx_i^{}]) \ge n_i^{}$ follow from the fact
that $\p \cap \KK[\xx_i^{}]$ contains the ideal of the projection of $X$ onto $\PP^{n_i}$.
According to Lemma~\ref{noether}, we have the inequality
\[
\height(\p) \ge \sum_{1 \le i \le q} \height(\p \cap \KK[\xx_i^{}]).
\]
This is equivalent to the inequality
\[
\sum_{k + 1 \le i \le q} n_i^{} \ge \sum_{k + 1 \le i \le q} \height(\p \cap \KK[\xx_i^{}]).
\]
This proves the claim.
The claim implies that $U_{k + 1}^{}$ is not contained in $\p$.
The same is true for all associated primes of $J$.
Since $\KK$ is infinite, we can choose $u_{k + 1}^{} \in U_{k + 1}^0$
such that $u_{k + 1}^{}$ does not lie in any associated prime of $J$.
Thus, $u_{k + 1}$ is a non-zerodivisor for $\SSS / J$,
hence $\{ u_1^{}, \dots, u_{k + 1}^{} \}$ is regular relative to $\SSS / I(X)$.
\end{proof}

\begin{lemma}
\label{minimal}
Assume that $\KK$ is algebraically closed.
Let $X \subset \VV$ be a zero-dimensional subscheme.
Write $\red(X) = \{ P_1^{}, \dots, P_m^{} \}$.
Choose $p \in \{ 1, \dots, q \}$ and choose $u_i^{} \in U_i^0$ for $1 \le i \le p$.
We claim that the ideals $\p_k^{} = (U_1^{}, \dots, U_p^{}) + I(P_k^{})$
for $1 \le k \le m$ are the minimal prime ideals containing $(u_1^{}, \dots, u_p^{}) + I(X)$.
\end{lemma}

\begin{proof}
Recall Remark~\ref{U_bis}.
Notice that $\p_k^{} = (U_1^{}, \dots, U_p^{}) + (U_{k, p+1}^{}, \dots, U_{k, q}^{})$.
This is clearly a prime ideal.
Some of these ideals may coincide, however,
if $\p_k^{} \neq \p_l^{}$, then $\p_k^{} \nsubseteq \p_l^{}$ and $\p_l^{} \nsubseteq \p_k^{}$.
The lemma reduces to proving that
\[
\bigcap_{1 \le k \le m} \p_k^{} = \rad((u_1^{}, \dots, u_p^{}) + I(X)).
\]
The inclusion ``$\supset$'' is obvious,
so we focus on proving the reverse inclusion.
We denote by $\rid$ the ideal on the r.h.s.
Take $f \in \bigcap_{1 \le k \le m} \p_k^{}$ and write $f = g + h$, where
\[
g \in (U_1^{}, \dots, U_p^{}) \quad \text{and} \quad
h \in \KK[x_{ij}^{} \mid p + 1 \le i \le q, \ 0 \le j \le n_i^{}].
\]
According to Lemma~\ref{radical}, $(U_1^{}, \dots, U_p^{}) \subset \rid$,
hence $g \in \rid$.
By construction,
\[
h \in \bigcap_{1 \le k \le m} (U_{ki}^{} \mid p + 1 \le i \le q) \subset \bigcap_{1 \le k \le m} I(P_k^{}) = \rad(I(X)) \subset \rid.
\]
We conclude that $f$ lies in $\rid$.
\end{proof}

\begin{lemma}
\label{operations}
Let $S$ be a $\ZZ$-graded $\KK$-algebra
and let $M$ be a $\ZZ$-graded $S$-module.
Assume that $\{ v_1^{}, \dots, v_p^{} \} \subset S$ is an $M$-regular sequence
and that all $v_i^{}$ are homogeneous of the same degree.
Consider a non-singular matrix $G = (\kappa_{ij}^{})_{1 \le i, j \le p}^{}$ with entries in $\KK$.
We claim that $\{ \kappa_{i1}^{} v_1^{} + \dots + \kappa_{ip}^{} v_p^{} \mid 1 \le i \le p \}$
constitutes an $M$-regular sequence.
\end{lemma}

\begin{proof}
If $G$ is a lower-triangular matrix,
then the claim simply follows from the definition of an $M$-regular sequence.
According to \cite[Corollary 17.5 and Theorem 17.6]{eisenbud_commutative},
a sequence $\{ w_1^{}, \dots, w_p^{} \} \subset S$ of homogeneous elements is $M$-regular
if and only if the Koszul complex $K(w_1^{}, \dots, w_p^{}) \tensor M$ is exact
(see Notation~\ref{complex}).
Permuting $\{ w_1^{}, \dots, w_p^{} \}$ results in isomorphic Koszul complexes.
Thus, every permutation of $\{ v_1^{}, \dots, v_p^{} \}$ remains an $M$-regular sequence.
The lemma follows from the fact
that the lower triangular matrices and the row permutations generate $\GL_p^{}(\KK)$.
\end{proof}

\begin{proposition}
\label{regular_2}
Assume that $\KK$ is algebraically closed.
Assume that $q \ge 3$.
Let $X \subset \VV$ be a zero-dimensional sub-ACM subscheme.
For $1 \le i \le q$ choose $u_i^{} \in U_i^0$.
We claim that there is $p \in \{ 2, \dots, q \}$
and there are $\kappa_i^{} \in \KK$ for $i \in \{ 2, \dots, q \} \setminus \{ p \}$
such that $\{ u_1^{} \} \cup \{ u_i^{} + \kappa_i^{} u_p^{} \mid i \in \{ 2, \dots, q \} \setminus \{ p \} \}$
is a regular sequence for $\SSS / I(X)$.
\end{proposition}

\begin{proof}
Write $\red(X) = \{ P_1^{}, \dots, P_m^{} \}$.
According to Lemma~\ref{minimal}, $\p_k^{} = (U_1^{}) + I(P_k^{})$ for $1 \le k \le m$
are the minimal prime ideals containing $(u_1^{}) + I(X)$.
Put $U' = \room \{ u_2^{}, \dots, u_q^{} \}$.
Performing induction on $l \in \{ 2, \dots, q - 1 \}$, we will construct a regular sequence
$\{ u_1^{}, v_2^{}, \dots, v_l^{} \}$ for $\SSS / I(X)$ with $v_i^{} \in U'$.
To start the induction, we consider the set $\A_1^{}$ of associated primes to $(u_1^{}) + I(X)$.
We claim that $U'$ is not contained in any $\p \in \A_1^{}$.
To prove this, we argue by contradiction.
Assume that $U' \subset \p$ and that $\p \in \A_1^{}$.
This ideal must contain one of the minimal associated primes to $(u_1^{}) + I(X)$,
say $\p_k^{} \subset \p$.
Thus, $U_1^{} \subset \p$ and $U_{ki}^{} \subset \p$ for $2 \le i \le q$.
It follows that
\begin{align*}
U & = U_1^{} + U_2^{} + \dots + U_q^{} \\
& = U_1^{} + \room \{ u_2^{}, U_{k2}^{} \} + \dots + \room \{ u_q^{}, U_{kq}^{} \} \\
& = U_1^{} + U_{k2}^{} + \dots + U_{kq}^{} + U' \subset \p,
\end{align*}
hence $\m = (U) \subset \p$,
so every element of $\m$ is a zerodivisor for $\SSS / ((u_1^{}) + I(X))$.
On the other hand, by Remark~\ref{U_0}, $u_1^{}$ is a non-zerodivisor for $\SSS / I(X)$, hence
\[
\depth(\SSS / ((u_1^{}) + I(X))) = \depth(\SSS / I(X)) - 1 = q - 2 \ge 1.
\]
We have reached a contradiction, which proves the claim.
We obtain a regular sequence $\{ u_1^{}, v_2^{} \}$ relative to $\SSS / I(X)$
by choosing $v_2^{} \in U' \setminus \bigcup_{\p \in \A_1} (\p \cap U')$.

We now perform the induction step.
Assume that $l \in \{ 2, \dots, q - 2 \}$ and that $\{ u_1^{}, v_2^{}, \dots, v_l^{} \}$ has already been constructed.
We denote by $\A_l^{}$ the set of associated primes to $(u_1^{}, v_2^{}, \dots, v_l^{}) + I(X)$.
Arguing as above, we can prove that $U'$ is not contained in any $\p$ from $\A_l^{}$.
Indeed, $\p_k^{} \subset \p$ for some $k$,
so, if $U' \subset \p$, then $U \subset \p$.
It would follow that every element of $\m$ is a zerodivisor for $\SSS / ((u_1^{}, v_2^{}, \dots, v_l^{}) + I(X))$.
On the other hand, this ring has depth $q - 1 - l \ge 1$.
This would yield a contradiction.
Choosing $v_{l + 1}^{} \in U' \setminus \bigcup_{\p \in \A_l} (\p \cap U')$
we obtain a regular sequence $\{ u_1^{}, v_2^{}, \dots, v_{l + 1}^{} \}$ relative to $\SSS / I(X)$.
This completes the induction step.

Thus far we have constructed a regular sequence $\{ u_1^{}, v_2^{}, \dots, v_{q - 1}^{} \}$ relative to $\SSS / I(X)$
such that $v_2^{}, \dots, v_{q - 1}^{}$ are linearly independent vectors in $U'$.
We write $v_l^{} = \sum_{2 \le i \le q} \lambda_{li}^{} u_i^{}$.
The matrix $\Lambda = (\lambda_{li}^{})_{2 \le l \le q - 1, \, 2 \le i \le q}^{}$ has maximal rank.
To simplify notations, we assume that the minor obtained by deleting the last column of $\Lambda$ is non-zero.
We now apply Lemma~\ref{operations} to the $\ZZ$-graded ring $\SSS$,
to the $\ZZ$-graded module $M = \SSS / ((u_1^{}) + I(X))$
and to the $M$-regular sequence $\{ v_2^{}, \dots, v_{q - 1}^{} \}$.
We take $G$ to be the inverse of $(\lambda_{li}^{})_{2 \le l, \, i \le q - 1}^{}$.
We obtain an $M$-regular sequence of the form $\{ u_i^{} + \kappa_i^{} u_q^{} \mid 2 \le i \le q - 1 \}$.
In general, if the minor obtained by deleting column $p$ of $\Lambda$ is non-zero,
then we obtain a regular sequence as in the proposition.
\end{proof}

\begin{proposition}
\label{regular_3}
Assume that $\KK$ is algebraically closed.
Let $X \subset \VV$ be a zero-dimensional ACM subscheme.
For $1 \le i \le q$ choose $u_i^{} \in U_i^0$.
We claim that $\{ u_1^{}, \dots, u_q^{} \}$ is regular relative to $\SSS / I(X)$.
\end{proposition}

\begin{proof}
Performing induction on $i \in \{ 1, \dots, q \}$,
we will show that $\{ u_1^{}, \dots, u_i^{} \}$ is regular relative to $\SSS / I(X)$.
By Remark~\ref{U_0}, $u_1^{}$ is a non-zerodivisor for $\SSS / I(X)$.
Assume now that $i \in \{1, \dots, q -1 \}$ and that $\{ u_1^{}, \dots, u_i^{} \}$ is regular for $\SSS / I(X)$.
By hypothesis, $\SSS / I(X)$ is Cohen-Macaulay,
hence $\SSS / ((u_1^{}, \dots, u_i^{}) + I(X))$ is Cohen-Macaulay,
and hence this ring is unmixed (see \cite[Corollary 18.14]{eisenbud_commutative}).
Thus, the associated primes of $(u_1^{}, \dots, u_i^{}) + I(X)$ are precisely the minimal primes of this ideal.
According to Lemma~\ref{minimal}, they are of the form
$\p_k^{} = (U_1^{}, \dots, U_i^{}) + I(P_k^{})$.
By construction, $u_{i + 1}^{}$ lies outside all ideals $\p_k^{}$,
hence $u_{i + 1}^{}$ is a non-zerodivisor of $\SSS / ((u_1^{}, \dots, u_i^{}) + I(X))$,
and hence $\{ u_1^{}, \dots, u_{i + 1}^{} \}$ is regular for $\SSS / I(X)$.
\end{proof}

\begin{lemma}
\label{intersection}
Let $S$ be a commutative ring and let $I \subset S$ be an ideal.
Assume that the sequence $\{ u_1^{}, \dots, u_p^{} \} \subset S$ is regular with respect to $S / I$.
We claim that
\[
(u_1^{}, \dots, u_p^{}) I = (u_1^{}, \dots, u_p^{}) \cap I.
\]
\end{lemma}

\begin{proof}
The inclusion ``$\subset$'' is obvious,
so we concentrate on proving the reverse inclusion.
We perform induction on $p$.
Assume that $p = 1$.
Take $f \in (u_1^{}) \cap I$
and write $f = u_1^{} g$.
In $S / I$ we have the relations $u_1^{} \hat{g} = \hat{u}_1^{} \hat{g} = \hat{f} = 0$.
By hypothesis, $u_1^{}$ is a non-zerodivisor for $S / I$,
hence $\hat{g} = 0$, that is, $g \in I$, and hence $f \in (u_1^{}) I$.
Assume that $p > 1$
and that the lemma is true for the $S / I$-regular sequence $\{ u_1^{}, \dots, u_{p - 1}^{} \}$.
Take $f \in (u_1^{}, \dots, u_p^{}) \cap I$
and write $f = u_1^{} g_1^{} + \dots + u_p^{} g_p^{}$.
In $S / ((u_1^{} + \dots + u_{p - 1}^{}) + I)$ we have the relations
\[
u_p^{} \hat{g}_p^{}
= \hat{u}_p^{} \hat{g}_p^{}
= \hat{f} - \hat{u}_1^{} \hat{g}_1^{} - \dots - \hat{u}_{p - 1}^{} \hat{g}_{p - 1}^{} = 0.
\]
By hypothesis, $u_p^{}$ is a non-zerodivisor for $S / ((u_1^{} + \dots + u_{p - 1}^{}) + I)$,
hence $\hat{g}_p^{} = 0$.
Write $g_p^{} = u_1^{} h_1^{} + \dots + u_{p - 1}^{} h_{p - 1}^{} + h_p^{}$, where $h_p^{} \in I$.
From the relation
\[
f - u_p^{} h_p^{} = u_1^{} (g_1^{} + u_p^{} h_1^{}) + \dots + u_{p - 1}^{} (g_{p - 1}^{} + u_p^{} h_{p - 1}^{})
\]
we see that $f - u_p^{} h_p^{} \in (u_1^{}, \dots, u_{p - 1}^{}) \cap I$.
By the induction hypothesis, this ideal is $(u_1^{}, \dots, u_{p - 1}) I$.
We conclude that $f \in (u_1^{}, \dots, u_p^{}) I$.
\end{proof}

\begin{lemma}
\label{regular_thrice}
Let $S$ be a commutative ring and let $I \subset S$ be an ideal.
Assume that the sequence $\{ u_1^{}, \dots, u_{p + 1}^{} \} \subset S$ is $S$-regular
and that the sequence $\{ u_1^{}, \dots, u_p^{} \}$ is $S / I$-regular.
We claim that $\{ u_1^{}, \dots, u_{p + 1}^{} \}$ is $I$-regular.
\end{lemma}

\begin{proof}
By hypothesis, $u_1^{}$ is a non-zerodivisor in $S$,
hence $u_1^{}$ is a non-zerodivisor for $I$.
Take $i \in \{ 1, \dots, p \}$.
We apply Lemma~\ref{intersection} to the $S / I$-regular sequence $\{ u_1^{}, \dots, u_i^{} \}$.
We deduce that $I / (u_1^{}, \dots, u_i^{})I$ is isomorphic, as an $S$-module, to an ideal of $S / (u_1^{}, \dots, u_i^{})$.
By hypothesis, $u_{i + 1}^{}$ is a non-zerodivisor for $S / (u_1^{}, \dots, u_i^{})$,
hence $u_{i + 1}^{}$ is a non-zerodivisor for $I / (u_1^{}, \dots, u_i^{})I$.
\end{proof}

\begin{proposition}
\label{regular_4}
Assume that $\KK$ is algebraically closed.
Let $X \subset \VV$ be a zero-dimensional sub-ACM subscheme.
For $1 \le i \le q$ choose $u_i^{} \in U_i^0$.
We claim that $\{ u_1^{}, \dots, u_q^{} \}$ is $I(X)$-regular.
\end{proposition}

\begin{proof}
Assume that $q = 2$.
According to Remark~\ref{U_0}, $\{ u_1^{} \}$ is regular for $\SSS / I(X)$.
Clearly, $\{ u_1^{}, u_2^{} \}$ is $\SSS$-regular.
From Lemma~\ref{regular_thrice} we deduce that $\{ u_1^{}, u_2^{} \}$ is $I(X)$-regular.
Assume that $q \ge 3$.
As in Proposition~\ref{regular_2}, let
\[
\{ w_1^{}, \dots, w_{q - 1}^{} \} = \{ u_i^{} + \kappa_i^{} u_p^{} \mid i \in \{ 1, \dots, q \} \setminus \{ p \} \}
\]
be a regular sequence relative to $\SSS / I(X)$.
Here $\kappa_1^{} = 0$.
Put $w_q^{} = u_p^{}$.
Clearly, $\{ w_1^{}, \dots, w_q^{} \}$ is $\SSS$-regular.
From Lemma~\ref{regular_thrice} we deduce that $\{ w_1^{}, \dots, w_q^{} \}$ is $I(X)$-regular.
Consider the column vectors
$\uu = (u_1^{}, \dots, u_q^{})^\trans$ and $\ww = (\ww_1^{}, \dots, \ww_q^{})^\trans$.
By construction, $\ww = \Lambda \uu$ for some $\Lambda \in \GL_q^{}(\KK)$.
We apply Lemma~\ref{operations} to the $I(X)$-regular sequence $\{ w_1^{}, \dots, w_q^{} \}$.
We take $G = \Lambda^{-1}$.
We conclude that $\{ u_1^{}, \dots, u_q^{} \}$ is $I(X)$-regular.
\end{proof}

\begin{proposition}
\label{regular_5}
Assume that $\KK$ is algebraically closed.
Let $X \subset \VV$ be a zero-dimensional ACM subscheme.
For $1 \le i \le q$ choose $u_i^{} \in U_i^0$.
For $1 \le i \le q$ choose $v_i^{} \in U_i^{} \setminus \KK u_i^{}$.
We claim that $\{ u_1^{}, \dots, u_q^{}, v_i^{} \}$ is $I(X)$-regular for every index $i$.
\end{proposition}

\begin{proof}
According to Proposition~\ref{regular_3}, $\{ u_1^{}, \dots, u_q^{} \}$ is regular relative to $\SSS / I(X)$.
Clearly, $\{ u_1^{}, \dots, u_q^{}, v_i^{} \}$ is $\SSS$-regular.
The proposition follows from Lemma~\ref{regular_thrice}.
\end{proof}

\begin{proposition}
\label{regular_6}
Assume that $\KK$ is infinite.
Assume that $X \subset \VV$ is a zero-dimensional ACM subscheme.
We claim that there are $u_i^{} \in U_i^{}$ and $v_i^{} \in U_i^{} \setminus \KK u_i^{}$
such that $\{ u_1^{}, \dots, u_q^{}, v_i^{} \}$ is $I(X)$-regular for every index $i$.
\end{proposition}

\begin{proof}
Proposition~\ref{regular_1} provides an $\SSS / I(X)$-regular sequence $\{ u_1^{}, \dots, u_q^{} \}$.
Clearly, $\{ u_1^{}, \dots, u_q^{}, v_i^{} \}$ is $\SSS$-regular.
The proposition follows from Lemma~\ref{regular_thrice}.
\end{proof}

\begin{proposition}
\label{regular_7}
Assume that $\KK$ is infinite.
Let $X \subset \VV$ be a zero-dimensional subscheme.
For $1 \le i \le q$ choose $u_i^{} \in U_i^0$ and $v_i^{} \in U_i^{} \setminus \KK u_i^{}$.
We claim that $\{ u_i^{}, v_j^{} \}$ is $I(X)$-regular for all indices $i$ and $j$.
\end{proposition}

\begin{proof}
According to Remark~\ref{U_0}, $\{ u_i^{} \}$ is regular for $\SSS / I(X)$.
Clearly, $\{ u_i^{}, v_j^{} \}$ is $\SSS$-regular.
The proposition follows from Lemma~\ref{regular_thrice}.
\end{proof}


\section{Finite ACM schemes}
\label{ACM}

\noindent
In this section we assume that the ground field $\KK$ is infinite.
We recall from section~\ref{regular} the notion of an ACM zero-dimensional scheme.
Lemma~\ref{concentrated} provides a class of examples of ACM schemes.
Recall that all zero-dimensional subschemes $X \subset \PP^n$ are ACM.
The next simplest case, when $\VV = \PP^1 \times \PP^1$,
was investigated by Giuffrida et al.\ \cite{giuffrida_postulation}.
For later use, we cite below some of the results in \cite[Section 4]{giuffrida_postulation}.

\begin{notation}
\label{characteristic}
Let $T \subset \ZZ^q$ be a subset.
The characteristic function $\X_T^{} \colon \ZZ^q \to \ZZ$ is given by
\[
\X_T^{}(a) = 
\begin{cases}
1 & \text{if $a \in T$}, \\
0 & \text{if $a \in \ZZ^q \setminus T$}.
\end{cases}
\]
\end{notation}

\begin{theorem}[Giuffrida et al.]
\label{giuffrida}
Consider the biprojective space $\VV = \PP^1 \times \PP^1$ over $\CC$,
with $\ZZ^2$-graded coordinate ring $\CC[x_0^{}, x_1^{}, y_0^{}, y_1^{}]$.
Let $X \subset \VV$ be a zero-dimensional subscheme.
Recall the relevant domain $R(X) = [0, s_1^{} - 1] \times [0, s_2^{} - 1]$ from Proposition~\ref{regularity_3}(v).
We claim that the following statements are equivalent:
\begin{enumerate}
\item[\emph{(i)}]
$X$ is ACM;
\item[\emph{(ii)}]
there is an integer $m \ge 0$
and there are $c_1^{}, \dots, c_m^{} \in R(X)$
such that the quasi-rectangular domain
\[
Q(X) = R(X) \setminus \bigcup_{1 \le k \le m} [c_k^{}, (s_1^{} - 1, s_2^{} - 1)] \subset \ZZ^2
\]
satisfies the condition $\der \hilb_X^{} = \X_{Q(X)}^{}$ (see Notation~\ref{characteristic});
\item[\emph{(iii)}]
there is an integer $m \ge 0$
and there are homogeneous polynomials $u_1^{}, \dots, u_{m + 1}^{} \in \CC[x_0^{}, x_1^{}]$
and $v_1^{}, \dots, v_{m + 1}^{} \in \CC[y_0^{}, y_1^{}]$
such that $I(X) = (v_1^{} \cdots v_k^{} u_{k + 1}^{} \cdots u_{m + 1}^{} \mid 0 \le k \le m + 1)$.
\end{enumerate}
Under the above conditions, $\deg(v_1^{} \cdots v_k^{} u_{k + 1}^{} \cdots u_{m + 1}^{}) = c_k^{}$.
Here $c_0^{} = (s_1^{}, 0)$ and $c_{m + 1}^{} = (0, s_2^{})$.
\end{theorem}

\noindent
As per Definition~\ref{relevance}, $Q(X)$ is a relevant quasi-rectangular domain for $\hilb_X^{}$,
in fact, the smallest possible quasi-rectangular relevant domain.

Other characterizations of the ACM property for zero-dimensional subschemes $X \subset \PP^1 \times \PP^1$ are known
in the case when $X$ is reduced,
see \cite[Theorem 4.8]{tuyl_hilbert}, \cite[Corollary 7.5]{marino_conductor}, \cite[Theorem 6.7]{marino_quadric},
\cite[Theorem 4.3]{guardo_points} and \cite[Theorem 8]{guardo_classifying},
and in the case when $X$ is a union of fat points,
see \cite[Theorem 2.1]{guardo_quadric} and \cite[Theorem 4.8]{guardo_hilbert}.
For other ambient spaces the focus has been entirely on reduced ACM schemes,
see \cite[Theorem 4.5 and Theorem 5.7]{guardo_points}, \cite[Theorem 3.16]{favacchio_points},
\cite[Proposition 3.2 and Theorem 3.7]{favacchio_multiprojective}.

We consider the algebra $\CC[x_0^{}, x_1^{}, y_0^{}, y_1^{}] / (x_0^{}, y_0^{}) = \CC[x_1^{}, y_1^{}]$
and we equip it with the inherited $\ZZ^2$-grading.
To wit, $\deg(x_1^{}) = e_1^{}$ and $\deg(y_1^{}) = e_2^{}$.
Condition (ii) from Theorem~\ref{giuffrida} is equivalent to saying
that $\der \hilb_X^{}$ is the Hilbert function of an artinian quotient of $\CC[x_1^{}, y_1^{}]$ by a monomial ideal,
i.e.\ an artinian $\ZZ^2$-graded quotient of $\CC[x_1^{}, y_1^{}]$.
This statement was partially generalized by Van~Tuyl in \cite{tuyl_hilbert}.
We consider the algebra
\[
\SSS_0^{} = \SSS / (x_{10}^{}, \dots, x_{q0}^{}) = \KK[x_{ij}^{} \mid 1 \le i \le q, \ 1 \le j \le n_i^{}]
\]
and we equip it with the induced $\ZZ^q$-grading.
Specifically, $\deg(x_{ij}^{}) = e_i^{}$.
According to \cite[Theorem 3.11]{tuyl_hilbert}, if $\KK$ is algebraically closed and
if $X \subset \VV$ is a zero-dimensional reduced and ACM subscheme,
then $\der \hilb_X^{}$ is the Hilbert function of an artinian $\ZZ^q$-graded quotient of $\SSS_0^{}$.
Conversely, for any artinian $\ZZ^q$-graded quotient $A$ of $\SSS_0^{}$,
there exists a zero-dimensional reduced ACM subscheme $X \subset \VV$
such that $\der \hilb_X^{} = \hilb_A^{}$.
The purpose of this section is to provide a version of this result
that does not require $X$ to be reduced or $\KK$ to be algebraically closed.

\begin{lemma}
\label{concentrated}
Assume that the subscheme $X \subset \VV$ is concentrated at a point
and that $I(X)$ is a monomial ideal.
We claim that $X$ is ACM.
\end{lemma}

\begin{proof}
By hypothesis, $\red(X) = \{ P \}$ for a closed point $P \in \VV$.
By the multiprojective version of Hilbert's Nullstellensatz, $I(P) = \rad(I(X))$.
As the radical of a monomial ideal, $I(P)$ itself is monomial.
But $I(P)$ is also a prime ideal,
hence $I(P)$ is generated by a subset of the set of variables.
We may assume that $I(P) = (x_{ij}^{} \mid 1 \le i \le q,\ 1 \le j \le n_i^{})$.
If $x_{i0}^{} u \in I(X)$ for a monomial $u$,
then, since $x_{i0}^{}$ does not vanish at $P$,
$u$ must lie in $I(X)$.
This shows that the minimal generators of $I(X)$ are monomials in the same variables that generate $I(P)$.
It has now become clear that $\{ x_{i0}^{} \mid 1 \le i \le q \}$ constitutes a regular sequence for $\SSS / I(X)$,
hence $\depth(\SSS / I(X)) \ge q$, and hence $\depth(\SSS / I(X)) = q$.
\end{proof}

\noindent
In the sequel we shall use Macaulay's theorem.
This theorem is usually stated for homogeneous ideals of polynomial rings,
see for instance \cite[Theorem 15.3]{eisenbud_commutative},
but it can easily be extended to the $\ZZ^q$-graded setting.

\begin{notation}
\label{leader}
Let us fix a monomial well-ordering on $\SSS$.
This is a well-ordering ``$\le$'' on the set $\MM$ of monic monomials of $\SSS$
which is compatible with multiplication: if $u, v, w \in \MM$ and $u < v$, then $u w < v w$.
We say that a monomial $u \in \MM$ occurs in a polynomial $f \in \SSS$
if $\kappa u$ is one of the monomials of $f$ for some $\kappa \in \KK \setminus \{ 0 \}$.
For $f \in \SSS \setminus \{ 0 \}$ we denote by $\lead(f)$ the largest $u \in \MM$ that occurs in $f$.
For an ideal $I \subset \SSS$ we introduce the \emph{leading ideal} $\lead(I) = ( \lead(f) \mid f \in I)$.
\end{notation}

\begin{theorem}[Macaulay]
\label{macaulay}
Let $I \subset \SSS$ be a $\ZZ^q$-homogeneous ideal.
Choose a monomial well-ordering on $\SSS$.
We claim that $\lead(I)$ is $\ZZ^q$-homogeneous
and that it has the same Hilbert function as $I$.
\end{theorem}

\begin{theorem}
\label{artinian}
Assume that $\KK$ is infinite.
Let $X \subset \VV$ be a zero-dimensional ACM subscheme.
We claim that $\der \hilb_X^{} = \hilb_A^{}$ for an artinian $\ZZ^q$-graded quotient $A$ of $\SSS_0^{}$.
Conversely, for any artinian $\ZZ^q$-graded quotient $A$ of $\SSS_0^{}$, we claim
that there exists a zero-dimensional ACM subscheme $X \subset \VV$
such that $\der \hilb_X^{} = \hilb_A^{}$.
\end{theorem}

\begin{proof}
Let $X \subset \VV$ be a zero-dimensional ACM subscheme.
By virtue of Proposition~\ref{regular_1},
there exists a regular sequence $\{ u_1^{}, \dots, u_q^{} \}$ for $\SSS / I(X)$
with $u_i^{} \in U_i^{}$.
In view of Lemma~\ref{koszul},
$\der \hilb_X^{}$ is the Hilbert function of $A = \SSS / ((u_1^{}, \dots, u_q^{}) + I(X))$.
Performing a linear change of coordinates on each $\PP^{n_i}$,
we may assume that $u_i^{} = x_{i0}^{}$ for all indices $i$,
so $A$ can be regarded as a $\ZZ^q$-graded quotient of $\SSS_0^{}$.
According to Theorem~\ref{main}, $\hilb_A^{}$ vanishes outside a rectangular domain,
hence $\dim_\KK^{} A$ is finite,
and hence $A$ is an artinian $\KK$-algebra.

Conversely, assume we are given an artinian $\ZZ^q$-graded algebra $A = \SSS_0^{} / I_0^{}$.
We choose a monomial well-ordering on $\SSS_0^{}$ (see Notation~\ref{leader})
and we apply Theorem~\ref{macaulay} to the $\ZZ^q$-homogeneous ideal $I_0^{}$.
We find a monomial ideal $J_0^{} = \lead(I_0^{})$ (see Notation~\ref{leader})
such that $\hilb_A^{} = \hilb_{\SSS_0 / J_0}$.
In particular, $\SSS_0^{} / J_0^{}$ is artinian,
hence $\rad(J_0^{}) = (x_{ij}^{} \mid 1 \le i \le q, \, 1 \le j \le n_i^{})$.
Let $J \subset \SSS$ be the ideal generated by $J_0^{}$.
Since $J$ is generated by monomials that do not involve the variables $x_{i0}$, $1 \le i \le q$,
it is obvious that $J$ is saturated.
Thus, $J = I(X)$ for a zero-dimensional subscheme $X \subset \VV$
which is concentrated on the point given by the ideal $(x_{ij}^{} \mid 1 \le i \le q, \, 1 \le j \le n_i^{})$.
According to Lemma~\ref{concentrated}, $X$ is ACM.
Since $J$ is generated by monomials that do not involve the variables $x_{i0}$, $1 \le i \le q$,
it is obvious that $\{ x_{10}^{}, \dots, x_{q0}^{} \}$ constitutes a regular sequence for $\SSS / J$.
Applying Lemma~\ref{koszul}, we find that $\der \hilb_X^{}$ is the Hilbert function of
$\SSS / ((x_{10}^{}, \dots, x_{q0}^{}) + J) = \SSS_0^{} / J_0^{}$.
\end{proof}

\noindent
The first claim of the above theorem, in the particular case when $\KK$ is algebraically closed and $X$ is reduced,
was obtained by Van~Tuyl.
Consult \cite[Theorem 3.11]{tuyl_hilbert}.
The second claim of the above theorem, in the particular case when $\KK$ is algebraically closed,
follows from Van~Tuyl's result.
Indeed, he proved that for any $A$ we can find a zero-dimensional reduced ACM subscheme $X \subset \VV$
such that $\der \hilb_X^{} = \hilb_A^{}$.

\begin{corollary}
\label{corners}
Assume that $\KK$ is infinite.
Let $X \subset (\PP^1)^q$ be a zero-dimensional ACM subscheme.
We claim that there exists a quasi-rectangular domain $Q(X) \subset \ZZ^q$
such that $\der \hilb_X^{} = \X_{Q(X)}^{}$ (see Notation~\ref{characteristic}).

Conversely, we claim that for any quasi-rectangular domain $Q \subset \ZZ^q$
there exists a zero-dimensional ACM subscheme $X \subset (\PP^1)^q$
such that $\der \hilb_X^{} = \X_Q^{}$.
\end{corollary}

\begin{proof}
We have $n_i^{} = 1$ for all indices $i$,
hence $\SSS_0^{} = \KK [x_{i1}^{} \mid 1 \le i \le q]$ with $\deg(x_{i1}^{}) = e_i^{}$.
An ideal of $\SSS_0^{}$ is $\ZZ^q$-homogeneous if and only if it is monomial.
Consider an artinian quotient $A = \SSS_0^{} / I_0^{}$ by a monomial ideal.
Since $A$ is artinian, $I_0^{}$ contains minimal generators of the form $x_{11}^{s_1}, \dots, x_{q1}^{s_q}$.
Let $c_1^{}, \dots, c_m^{}$ be the degrees of the remaining minimal generators of $I_0^{}$, if any.
We have $\hilb_A^{} = \X_Q^{}$, where
\[
Q = [0, s_1^{} - 1] \times \dots \times [0, s_q^{} - 1]
\setminus \bigcup_{1 \le k \le m} [c_k^{}, (s_1^{} - 1, \dots, s_q^{} - 1)].
\]
Conversely, for any quasi-rectangular domain $Q \subset \ZZ^q$,
we can find an artinian $\ZZ^q$-graded quotient $A = \SSS_0^{} / I_0^{}$
such that $\hilb_A^{} = \X_Q^{}$.
\end{proof}

\noindent
The first claim of the above corollary is a generalization
of the implication ``(i)$\Longrightarrow$(ii)'' of Theorem~\ref{giuffrida}.
The first claim of the above corollary, in the particular case when $\KK$ is algebraically closed and $X$ is reduced,
was obtained by Van~Tuyl.
Consult \cite[Corollary 3.14]{tuyl_hilbert}.
The second claim of the above corollary, in the particular case when $\KK$ is algebraically closed,
follows from Van~Tuyl's result.
Indeed, he proved that for any $Q$ we can find a zero-dimensional reduced ACM subscheme $X \subset (\PP^1)^q$
such that $\der \hilb_X^{} = \X_Q^{}$.


\section{Further constraints on the Hilbert functions}
\label{constraints}

\noindent
In this section we assume that the ground field $\KK$ is infinite.
This section is devoted to a better understanding of the problem of classification
of the functions $\ZZ^q \to \ZZ$ that arise as Hilbert functions of zero-dimensional subschemes $X \subset \VV$.
A classical theorem of Macaulay (see \cite[Theorem 4.2.10]{bruns-herzog})
provides a classification of the Hilbert functions of $\ZZ$-graded $\KK$-algebras.
A recent theorem of Favacchio (see \cite[Theorem 4.8]{favacchio_numerical})
provides a classification of the Hilbert functions of $\ZZ^2$-graded $\KK$-algebras.
Invoking Theorem~\ref{artinian}, we obtain a characterization of the functions $\hilb_X^{}$,
for zero-dimensional subschemes $X \subset \PP^n$, respectively,
for zero-dimensional ACM subschemes $X \subset \PP^{n_1} \times \PP^{n_2}$.
The problem of describing the functions $\hilb_X^{}$ in the case when $X \subset \PP^{n_1} \times \PP^{n_2}$
is sub-ACM, or in the case when $q \ge 3$, remains open.
In this section we make progress on this problem by exhibiting certain conditions
that the functions $\hilb_X^{}$ and $\hilb_{I(X)}^{}$ must satisfy.
These constraints are formulated in terms of the partial difference functions, defined below.
The emphasis will be on ACM and sub-ACM schemes.
All constraints arise in the manner of Theorem~\ref{artinian}:
we exploit the regular sequences from section~\ref{regular},
and then we apply Lemma~\ref{koszul}.
At the end of the section we give a second proof to Theorem~\ref{main}
in the particular case when $\VV = (\PP^1)^q$ and $X$ is ACM or sub-ACM.

Let $\F \colon \ZZ^q \to \ZZ$ be a function.
For $1 \le i \le q$ we consider the \emph{partial difference} function
\[
\frac{\der \F}{\der a_i^{}} \colon \ZZ^q \lra \ZZ
\quad \text{given by} \quad
\frac{\der \F}{\der a_i^{}}(a) = \F(a) - \F(a - e_i^{}).
\]
We write
\[
\frac{\der}{\der a_{i_1}^{}} \cdots \frac{\der}{\der a_{i_p}^{}} \F
= \frac{\der^p \F}{\der a_{i_1}^{} \cdots \der a_{i_p}^{}}
\quad \text{and} \quad
\frac{\der^p \F}{\der a_i^{} \cdots \der a_i^{}}
= \frac{\der^p \F}{\der a_i^p}.
\]
Notice that
\[
\der \F = \frac{\der^q \F}{\der a_1^{} \cdots \der a_q^{}}
\quad \text{and} \quad
\frac{\der^p}{\der a^p} {a + n \choose n} = {a + n - p \choose n - p}.
\]
We will use the abbreviation
\[
\frac{\der \der \F}{\der a_i^{}} = \frac{\der^{q + 1} \F}{\der a_1^{} \cdots \der a_q^{} \der a_i^{}}.
\]
We have
\[
\hilb_\SSS^{}(a) = \prod_{1 \le i \le q} {a_i^{} + n_i^{} \choose n_i^{}}.
\]
We obtain the equation
\begin{equation}
\label{difference_1}
\der \hilb_\SSS^{}(a) = \prod_{1 \le i \le q} {a_i^{} + n_i^{} - 1 \choose n_i^{} - 1}.
\end{equation}
We have
\[
\frac{\der \der \hilb_\SSS^{}}{\der a_i^{}}(a)
= {a_i^{} + n_i^{} - 2 \choose n_i^{} - 2} \prod_{\substack{1 \le j \le q \\ j \neq i}} {a_j^{} + n_j^{} - 1 \choose n_j^{} - 1}.
\]

\begin{proposition}
\label{constraint_1}
Assume that $\KK$ is algebraically closed.
Assume that the zero-dimensional subscheme $X \subset \VV$ is sub-ACM.
We make the following claims:
\begin{enumerate}
\item[\emph{(i)}]
$\der \hilb_{I(X)}^{} \ge 0$;
\item[\emph{(ii)}]
$\der \hilb_{I(X)}^{}(a) > 0$ if $I(X)_a^{} \neq \{ 0 \}$;
\item[\emph{(iii)}]
$\der \hilb_X^{} \le \der \hilb_\SSS^{}$;
\item[\emph{(iv)}]
if $\der \hilb_X^{}(a) = \der \hilb_\SSS^{}(a)$ for some $a \in \ZZ_+^q$,
then $\hilb_X^{} = \hilb_\SSS^{}$ on $[0, a]$;
\item[\emph{(v)}]
$\dfrac{\der^p \hilb_{I(X)}^{}}{\der a_{i_1}^{} \cdots \der a_{i_p}^{}} \ge 0$ for all indices $1 \le i_1^{} < \dots < i_p^{} \le q$;
\item[\emph{(vi)}]
$\dfrac{\der^p \hilb_{I(X)}^{}}{\der a_{i_1}^{} \cdots \der a_{i_p}^{}}(a) > 0$
if $I(X)_a^{} \neq \{ 0 \}$ and $1 \le i_1^{} < \dots < i_p^{} \le q$.
\end{enumerate}
\end{proposition}

\begin{proof}
Recall Notation~\ref{U}.
As per Proposition~\ref{regular_4},
we can construct $I(X)$-regular sequences $\{ u_1^{}, \dots, u_q^{} \}$ with generic $u_i^{} \in U_i^{}$.
Write $N = I(X) / (u_1^{}, \dots, u_q^{}) I(X)$.
Applying Lemma~\ref{koszul}, we deduce that $\der \hilb_{I(X)}^{} = \hilb_N^{}$.
This function takes only non-negative values.
We have proved claim (i).
The same proof applies to claim (v),
except that we consider the $I(X)$-regular sequence $\{ u_{i_1}^{}, \dots, u_{i_p}^{} \}$.

If $I(X)_a^{} \neq \{ 0 \}$, then we can choose $u_i^{}$
such that the subvariety given by the ideal $(u_1^{}, \dots, u_q^{})$ is not contained in the zero-set of $I(X)_a$.
Thus, $N_a^{} \neq \{ 0 \}$, hence $\der \hilb_{I(X)}^{}(a) > 0$.
This proves claim (ii).
The same proof applies in claim (vi),
except that we consider the ideal $(u_{i_1}^{}, \dots, u_{i_p}^{})$.
Claim (iii) follows from the equation $\der \hilb_X^{} = \der \hilb_\SSS^{} - \der \hilb_{I(X)}^{}$ and from claim (i).

Assume that $\der \hilb_X^{}(a) = \der \hilb_\SSS^{}(a)$ for some $a \in \ZZ_+^q$.
Thus, $\der \hilb_{I(X)}^{}(a) = 0$.
From claim (ii) we deduce that $I(X)_a^{} = \{ 0 \}$.
A fortiori, $\hilb_{I(X)}^{} = 0$ on $[0, a]$,
hence $\hilb_X^{} = \hilb_\SSS^{}$ on $[0, a]$.
This proves claim (iv).
\end{proof}

\begin{proposition}
\label{constraint_2}
Assume that $\KK$ is infinite.
Assume that the zero-dimensional subscheme $X \subset \VV$ is ACM.
As in Theorem~\ref{artinian}, let $A$ be an artinian algebra such that $\der \hilb_X^{} = \hilb_A^{}$.
We make the following claims:
\begin{enumerate}
\item[\emph{(i)}]
$\dfrac{\der \der \hilb_{I(X)}^{}}{\der a_i^{}} \ge 0$ for all $i \in \{ 1, \dots, q \}$;
\item[\emph{(ii)}]
$\dfrac{\der \der \hilb_{I(X)}^{}}{\der a_i^{}}(a) > 0$ if $\KK$ is algebraically closed, $I(X)_a^{} \neq \{ 0 \}$ and $n_i^{} \ge 2$;
\item[\emph{(iii)}]
$\dfrac{\der \hilb_A^{}}{\der a_i^{}} \le \dfrac{\der \der \hilb_\SSS^{}}{\der a_i^{}}$
for all $i \in \{ 1, \dots, q \}$ such that $n_i^{} \ge 2$;
\item[\emph{(iv)}]
if $n_i^{} \ge 2$, $\KK$ is algebraically closed and
$\dfrac{\der \hilb_A^{}}{\der a_i^{}}(a) = \dfrac{\der \der \hilb_\SSS^{}}{\der a_i^{}}(a)$ for some $a \in \ZZ_+^q$,
then $\hilb_X^{} = \hilb_\SSS^{}$ on $[0, a]$;
\item[\emph{(v)}]
$\der \hilb_{I(X)}^{} \ge 0$;
\item[\emph{(vi)}]
$\der \hilb_{I(X)}^{}(a) > 0$ if $\KK$ is algebraically closed and $I(X)_a^{} \neq \{ 0 \}$;
\item[\emph{(vii)}]
$\hilb_A^{} \le \der \hilb_\SSS^{}$;
\item[\emph{(viii)}]
if $\KK$ is algebraically closed and $\hilb_A^{}(a) = \der \hilb_\SSS^{}(a)$ for some $a \in \ZZ_+^q$,
then $\hilb_X^{} = \hilb_\SSS^{}$ on $[0, a]$;
\item[\emph{(ix)}]
$\dfrac{\der^p \hilb_{I(X)}^{}}{\der a_{i_1}^{} \cdots \der a_{i_p}^{}} \ge 0$ for all indices $1 \le i_1^{} < \dots < i_p^{} \le q$;
\item[\emph{(x)}]
$\dfrac{\der^p \hilb_{I(X)}^{}}{\der a_{i_1}^{} \cdots \der a_{i_p}^{}}(a) > 0$
if $\KK$ is algebraically closed, $I(X)_a^{} \neq \{ 0 \}$ and $1 \le i_1^{} < \dots < i_p^{} \le q$.
\end{enumerate}
\end{proposition}

\begin{proof}
Consider the $I(X)$-regular sequence $\{ u_1^{}, \dots, u_q^{}, v_i^{} \}$ from Proposition~\ref{regular_6}.
Write $N = I(X) / (u_1^{}, \dots, u_q^{}, v_i^{}) I(X)$.
By analogy with Lemma~\ref{koszul}, we can prove that
$\dfrac{\der \der \hilb_{I(X)}^{}}{\der a_i^{}} = \hilb_N^{}$.
This function takes only non-negative values.
We have proved claim (i).
Assume that $I(X)_a^{} \neq \{ 0 \}$ and $n_i^{} \ge 2$.
According to Proposition~\ref{regular_5}, $u_1^{}, \dots, u_q^{}$ and $v_i^{}$ can be chosen generically.
We choose them in such a way that the subvariety given by the ideal
$(u_1^{}, \dots, u_q^{}, v_i^{})$ is not contained in the zero-set of $I(X)_a^{}$.
Thus, $N_a^{} \neq \{ 0 \}$, hence $\hilb_N^{}(a) > 0$.
This proves claim (ii).
Claim (iii) follows from the equation
\[
\dfrac{\der \hilb_A^{}}{\der a_i^{}}
= \dfrac{\der \der \hilb_\SSS^{}}{\der a_i^{}} - \dfrac{\der \der \hilb_{I(X)}^{}}{\der a_i^{}}
\]
and from claim (i).
To prove the remaining claims we can argue as in the proof of Proposition~\ref{constraint_1}.
Note also that claim (v) follows from claim (i) and from the formula
\[
\der \hilb_{I(X)}^{} = \sum_{0 \le k \le a_i} \dfrac{\der \der \hilb_{I(X)}^{}}{\der a_i^{}}(a - k e_i^{}).
\]
Likewise, in the case when $n_i^{} \ge 2$, claim (vi) follows from claims (i) and (ii),
and claim (vii) follows from claim (iii)
\end{proof}

\begin{corollary}
\label{constraint_3}
Assume that $\KK$ is algebraically closed.
Assume that the zero-dimensional subscheme $X \subset (\PP^1)^q$ is ACM or sub-ACM.
We make the following claims:
\begin{enumerate}
\item[\emph{(i)}]
$\der \hilb_X^{}(a) \le 1$ for all $a \in \ZZ_+^q$;
\item[\emph{(ii)}]
if $\der \hilb_X^{}(a) = 1$ for some $a \in \ZZ_+^q$,
then $\der \hilb_X^{} = 1$ on $[0, a]$.
\end{enumerate}
\end{corollary}

\begin{proof}
Substituting $n_i^{} = 1$ into formula~\eqref{difference_1}, we obtain $\der \hilb_\SSS^{} (a) = 1$ for $a \in \ZZ_+^q$.
Substituting this expression into Proposition~\ref{constraint_1}(iii) and Proposition~\ref{constraint_2}(vii)
yields claim (i).
Substituting this expression into Proposition~\ref{constraint_1}(iv) and Proposition~\ref{constraint_2}(viii)
yields claim (ii).
\end{proof}

\noindent
In the case when $X$ is ACM, the above corollary also follows from Corollary~\ref{corners}.
The above result, in the particular case when $\VV = \PP^1 \times \PP^1$, was obtained by Giuffrida et al.
Consult \cite[Proposition 2.7]{giuffrida_postulation}.

\begin{proposition}
\label{constraint_4}
Assume that $\KK$ is infinite.
Let $X \subset \VV$ be a zero-dimensional subscheme.
We make the following claims:
\begin{enumerate}
\item[\emph{(i)}]
$\dfrac{\der^2 \hilb_{I(X)}^{}}{\der a_i^{} \der a_j^{}} \ge 0$ for all indices $1 \le i \le j \le q$;
\item[\emph{(ii)}]
$\dfrac{\der^2 \hilb_{I(X)}^{}}{\der a_i^{} \der a_j^{}}(a) > 0$ if $\KK$ is algebraically closed, $I(X)_a^{} \neq \{ 0 \}$
and $1 \le i < j \le q$;
\item[\emph{(iii)}]
$\dfrac{\der^2 \hilb_{I(X)}^{}}{\der a_i^2}(a) > 0$ if $\KK$ is algebraically closed, $I(X)_a^{} \neq \{ 0 \}$ and $n_i^{} \ge 2$;
\item[\emph{(iv)}]
$\dfrac{\der^2 \hilb_X^{}}{\der a_i^{} \der a_j^{}} \le \dfrac{\der^2 \hilb_\SSS^{}}{\der a_i^{} \der a_j^{}}$
for all indices $1 \le i < j \le q$;
\item[\emph{(v)}]
if $\KK$ is algebraically closed and
$\dfrac{\der^2 \hilb_X^{}}{\der a_i^{} \der a_j^{}}(a) = \dfrac{\der^2 \hilb_\SSS^{}}{\der a_i^{} \der a_j^{}}(a)$
for some $a \in \ZZ_+^q$ and indices $1 \le i < j \le q$,
then $\hilb_X^{} = \hilb_\SSS^{}$ on $[0, a]$;
\item[\emph{(vi)}]
if $n_i^{} \ge 2$, $\KK$ is algebraically closed and
$\dfrac{\der^2 \hilb_X^{}}{\der a_i^2}(a) = \dfrac{\der^2 \hilb_\SSS^{}}{\der a_i^2}(a)$ for some $a \in \ZZ_+^q$,
then $\hilb_X^{} = \hilb_\SSS^{}$ on $[0, a]$;
\item[\emph{(vii)}]
$\dfrac{\der \hilb_X^{}}{\der a_i^{}} \ge 0$ for all $i \in \{ 1, \dots, q \}$.
\end{enumerate}
\end{proposition}

\begin{proof}
We use the $I(X)$-regular sequence $\{ u_i^{}, v_j^{} \}$ provided by Proposition~\ref{regular_7}
and we repeat the arguments from the proof of Proposition~\ref{constraint_1}.
Claim (vii) follows from the fact that $\{ u_i^{} \}$ is regular for $\SSS / I(X)$,
see Remark~\ref{U_0}.
\end{proof}

\noindent
As an application of the above results, we will give a second proof to a particular case of Theorem~\ref{main}.
We formulate this as a separate proposition.

\begin{proposition}
\label{lines_1}
Assume that $\KK$ is algebraically closed.
Assume that the zero-dimensional subscheme $X \subset (\PP^1)^q$ is ACM or sub-ACM.
For each index $i \in \{ 1, \dots, q \}$,
let $s_i^{}$ be the length of the projection of $X$ onto the $i$-th copy of $\PP^1$.
We claim that $[0, s_1^{} - 1] \times \dots \times [0, s_q^{} - 1]$
is the smallest rectangular relevant domain for $\hilb_X^{}$.
\end{proposition}

\begin{proof}
In view of Lemma~\ref{relevant},
we must show that $\hilb_X^{}(a) = \hilb_X^{}(a - (a_i^{} - s_i^{} + 1) e_i^{})$ if $a_i^{} \ge s_i^{} - 1$.
Equivalently, we must show that $\dfrac{\der \hilb_X^{}}{\der a_i^{}}(a) = 0$ if $a_i^{} \ge s_i^{}$.
By symmetry, it is enough to prove the statement
\[
\frac{\der \hilb_X^{}}{\der a_1^{}}(a) = 0 \quad \text{if} \ a \in \ZZ_+^q \ \text{and} \ a_1^{} \ge s_1^{}.
\]
Given $a \in \ZZ_+^q$, we put $\sigma(a) = a_2^{} + \dots + a_q^{}$.
We perform induction on $\sigma(a)$.
To begin the induction, we assume that $\sigma(a) = 0$.
Write $Z_1^{} = \pr_1^{}(X)$.
According to Proposition~\ref{regularity_2}(i) and (v),
\[
\hilb_X^{} (b_1^{}, 0 , \dots, 0) = \hilb_{Z_1}^{}(b_1^{}) = s_1^{} \quad \text{for} \quad b_1^{} \ge s_1^{} - 1.
\]
This leads to the desired outcome:\ $\dfrac{\der \hilb_X^{}}{\der a_1^{}}(a) = 0$.
We now perform the induction step.
We assume that $\sigma(a) > 0$.
To simplify notations, we assume that $a_2^{}, \dots, a_p^{}$ are positive
and $a_{p + 1}^{}, \dots, a_q^{}$ are zero for some $p \in \{ 2, \dots, q \}$.
From the definition of the partial difference functions we easily obtain the formula
\[
\frac{\der \hilb_X^{}}{\der a_1^{}}(a)
= \frac{\der^p \hilb_X^{}}{\der a_1^{} \cdots \der a_p^{}}(a)
+ \sum_{1 \le k \le p - 1} \frac{\der^k \hilb_X^{}}{\der a_1^{} \cdots \der a_k^{}}(a - e_{k + 1}^{}).
\]
Since $\sigma(a - e_{k + 1}^{}) < \sigma(a)$, 
$\dfrac{\der^k \hilb_X^{}}{\der a_1^{} \cdots \der a_k^{}}(a - e_{k + 1}^{})$
is a finite sum of expressions of the form $\pm \dfrac{\der \hilb_X^{}}{\der a_1^{}}(b)$,
with $\sigma(b) < \sigma(a)$ and with $b_1^{} = a_1^{}$.
By the induction hypothesis these expressions vanish.
We are led to the equations
\begin{align*}
\frac{\der \hilb_X^{}}{\der a_1^{}}(a)
& = \frac{\der^p \hilb_X^{}}{\der a_1^{} \cdots \der a_p^{}}(a) \\
& = (a_{p + 1}^{} + 1) \cdots (a_q^{} + 1) - \frac{\der^p \hilb_{I(X)}^{}}{\der a_1^{} \cdots \der a_p^{}}(a) \\
& = 1 - \frac{\der^p \hilb_{I(X)}^{}}{\der a_1^{} \cdots \der a_p^{}}(a).
\end{align*}
We know that $I(Z_1^{})_{s_1}^{} \neq \{ 0 \}$.
It follows that $I(X)_a^{} \neq \{ 0 \}$.
Since we are assuming that $X$ is ACM or sub-ACM,
we may apply Proposition~\ref{constraint_1}(vi) and Proposition~\ref{constraint_2}(x)
in order to deduce that
\[
\frac{\der^p \hilb_{I(X)}^{}}{\der a_1^{} \cdots \der a_p^{}}(a) > 0.
\qquad \text{A fortiori}, \quad \frac{\der \hilb_X^{}}{\der a_1^{}}(a) \le 0.
\]
According to Proposition~\ref{constraint_4}(vii),
the reverse inequality $\dfrac{\der \hilb_X^{}}{\der a_1^{}}(a) \ge 0$ also holds.
We obtain the desired outcome:\ $\dfrac{\der \hilb_X^{}}{\der a_1^{}}(a) = 0$.
This concludes the proof of the induction step.
\end{proof}

\noindent
The above line of argument, in the particular case when $\VV = \PP^1 \times \PP^1$, is due to Giuffrida et al.
Consult \cite[Remark 2.8 and Theorem 2.11]{giuffrida_postulation}.
We have adapted their proof to the case of arbitrary $q$.
In the case when $q = 2$ there is no restriction on $X$ because
every zero-dimensional subscheme $X \subset \PP^1 \times \PP^1$ is ACM or sub-ACM.


\section{A vanishing result for $\der \hilb$}
\label{formula}

\noindent
In this section we assume that $\VV = (\PP^1)^q$.
We write $\SSS = \KK[x_1^{}, y_1^{}, \dots, x_q^{}, y_q^{}]$,
where $\deg(x_i^{}) = e_i^{}$ and $\deg(y_i^{}) = e_i^{}$.
We saw at Theorem~\ref{giuffrida} and at Corollary~\ref{corners}
that zero-dimensional ACM subschemes $X \subset \VV$ that are not complete intersections
have a quasi-rectangular relevant domain $Q(X)$ 
that is strictly contained in the rectangular relevant domain $R(X)$ introduced at Theorem~\ref{main}.
This section and the next are devoted to finding a procedure (Proposition~\ref{procedure})
for constructing a quasi-rectangular relevant domain $D(X) \subset R(X)$
that applies to schemes $X$ which are not necessarily ACM.
We restrict our attention only to schemes $X$ for which $I(X)$ is a monomial ideal.
The domain $D(X)$ may coincide with $R(X)$ or may be strictly contained in $R(X)$, depending on the scheme.
At the end of section~\ref{end} we shall give examples
in which $D(X)$ is strictly contained in $R(X)$.

In this section we do some preparatory work.
We obtain a vanishing criterion for $\der \hilb_{\SSS / J}^{}$,
where $J \subset \SSS$ is a monomial ideal.
In order to achieve this, we need to take two preliminary steps.
Firstly, at Lemma~\ref{delta} we obtain a combinatorial formula for $\der \hilb_{\SSS / J}^{}$.
This formula actually holds for any $\ZZ^q$-graded polynomial ring,
i.e.\ for arbitrary values of $n_1^{}, \dots, n_q^{}$.
The second step, Lemma~\ref{gammas}, is also combinatorial
and breaks down if there are more than two variables of degree $e_i^{}$.
This is the technical reason why our ambient space needs to be a product of projective lines.

We find it convenient to work with the function $\hilb_J^{}$.
Substituting $n_i^{} = 1$ into equation~\eqref{difference_1}
we get $\der \hilb_\SSS^{}(a) = 1$ for $a \in \ZZ_+^q$.
The formula $\hilb_{\SSS / J}^{} = \hilb_\SSS^{} - \hilb_J^{}$ yields
\begin{equation}
\label{difference_2}
\der \hilb_{\SSS / J}^{}(a) = 1 - \der \hilb_J^{}(a) \quad \text{for} \quad a \in \ZZ_+^q.
\end{equation}
Let $\MM$ be the set of monic monomials of $\SSS$.
Let $\Gamma(J) = \{ f_1^{}, \dots, f_m^{} \} \subset \MM$ be the set of minimal generators of $J$.
Fix $a \in \ZZ^q$.
For all integers $p \in \{ 1, \dots, m \}$ we write
\[
\Gamma_a^p(J) = \{ (f_{k_1}^{}, \dots, f_{k_p}^{}) \mid 1 \le k_1^{} < \dots < k_p^{} \le m, \
\deg(\lcm(f_{k_1}^{}, \dots, f_{k_p}^{})) \le a \}.
\]

\begin{lemma}
\label{delta}
We consider $a \in \ZZ_+^q$.
We adopt the above notations.
We claim that
\[
\der \hilb_{\SSS / J}^{}(a) = 1 + \sum_{p \ge 1} (-1)^p \left| \Gamma_a^p(J) \right|.
\]
\end{lemma}

\begin{proof}
Given indices $1 \le k_1^{} < \dots < k_p^{} \le m$ we write
\[
d_{k_1 \dots k_p}^{} = (d_{k_1 \dots k_p}^1, \dots, d_{k_1 \dots k_p}^q)
= \deg(\lcm(f_{k_1}^{}, \dots, f_{k_p}^{})).
\]
By definition, for $b \in \ZZ^q$,
\[
\hilb_J^{}(b) = \Big| \bigcup_{1 \le k \le m} \{ u \in \MM \mid f_k^{} \ \text{divides} \ u, \ \deg(u) = b \} \Big|.
\]
Applying the inclusion-exclusion principle, we obtain the formula
\[
\hilb_J^{}(b) = \sum_{1 \le p \le m} (-1)^{p+1} \sum_{1 \le k_1^{} < \dots < k_p^{} \le m}
\big| \{ u \in \MM \mid f_{k_1}^{}, \dots, f_{k_p}^{} \ \text{divide} \ u, \ \deg(u) = b \} \big|.
\]
Ignoring the empty sets on the r.h.s., we calculate:
\[
\hilb_J^{}(b) = \sum_{1 \le p \le m} (-1)^{p+1} \sum_{\substack{1 \le k_1^{} < \dots < k_p^{} \le m \\ d_{k_1 \dots k_p} \le b}}
(b_1^{} - d_{k_1 \dots k_p}^1 + 1) \cdots (b_q^{} - d_{k_1 \dots k_p}^q + 1).
\]
Assume now that $b = a + c$, where $c_i^{} \in \{ -1, 0 \}$ for all indices $i \in \{ 1, \dots, q \}$.
If $d \in \ZZ^q$ satisfies the conditions $d \le a$ and $d \nleq b$,
then there is an index $i$ such that $d_i^{} = a_i^{} = b_i^{} + 1$,
forcing the equation $(b_1^{} - d_1^{} + 1) \cdots (b_q^{} - d_q^{} + 1) = 0$.
Thus, on the r.h.s.\ of the above formula we can add all the terms for which
$d_{k_1 \dots k_p}^{} \le a$ but $d_{k_1 \dots k_p}^{} \nleq b$,
i.e.\ we can replace $b$ by $a$ under the summation sign:
\[
\hilb_J^{}(b) = \sum_{1 \le p \le m} (-1)^{p+1} \sum_{\substack{1 \le k_1^{} < \dots < k_p^{} \le m \\ d_{k_1 \dots k_p} \le a}}
(b_1^{} - d_{k_1 \dots k_p}^1 + 1) \cdots (b_q^{} - d_{k_1 \dots k_p}^q + 1).
\]
This equation holds for $b = a + c$ for all possible $c \in \{ -1, 0 \}^q$,
hence we may apply the $\der$ operator calculated at $a$ to both sides:
\begin{align*}
\der \hilb_J^{}(a) & = \sum_{1 \le p \le m} (-1)^{p+1} \!\!\!\!\!\!\!
\sum_{\substack{1 \le k_1^{} < \dots < k_p^{} \le m \\ d_{k_1 \dots k_p} \le a}}
\der \big((b_1^{} - d_{k_1 \dots k_p}^1 + 1) \cdots (b_q^{} - d_{k_1 \dots k_p}^q + 1)\big)|_{b = a} \\
& = \sum_{1 \le p \le m} (-1)^{p+1} \!\!\!\!\!\!\!
\sum_{\substack{1 \le k_1^{} < \dots < k_p^{} \le m \\ d_{k_1 \dots k_p} \le a}} 1
\quad = \sum_{1 \le p \le m} (-1)^{p+1} \left| \Gamma_a^p(J) \right|.
\end{align*}
To conclude the proof of the lemma we employ equation~\eqref{difference_2}.
\end{proof}

\begin{lemma}
\label{gammas}
Assume that $2 \le p \le m$.
Under the above notations, we claim that
\begin{multline*}
\Gamma_a^p(J) = \{ (f_{k_1}^{}, \dots, f_{k_p}^{}) \mid 1 \le k_1^{} < \dots < k_p^{} \le m, \
(f_{k_\mu}^{}, f_{k_\nu}^{}) \in \Gamma_a^2(J) \\
\text{for all indices} \ 1 \le \mu < \nu \le p \}.
\end{multline*}
\end{lemma}

\begin{proof}
Assume that $(f_{k_1}^{}, \dots, f_{k_p}^{})$ lies in $\Gamma_a^p(J)$.
Then, for all indices $1 \le \mu < \nu \le p$,
\[
\deg(\lcm(f_{k_\mu}^{}, f_{k_\nu}^{})) \le \deg(\lcm(f_{k_1}^{}, \dots, f_{k_p}^{})) \le a.
\]
This proves the inclusion ``$\subset$''.
Conversely, assume that $(f_{k_1}^{}, \dots, f_{k_p}^{})$ belongs to the set on the r.h.s.
For each index $\mu \in \{ 1, \dots, p \}$ write
\[
f_{k_\mu}^{} = x_1^{\alpha_{\mu 1}} y_1^{\beta_{\mu 1}} \cdots x_q^{\alpha_{\mu q}} y_q^{\beta_{\mu q}}.
\]
For all indices $i \in \{ 1, \dots, q \}$ put $\alpha_i^{} = \max_{1 \le \mu \le p} \alpha_{\mu i}^{}$
and $\beta_i^{} = \max_{1 \le \mu \le p} \beta_{\mu i}^{}$.
Thus,
\[
\lcm(f_{k_1}^{}, \dots, f_{k_p}^{}) = x_1^{\alpha_1} y_1^{\beta_1} \cdots x_q^{\alpha_q} y_q^{\beta_q}.
\]
For a fixed index $i \in \{ 1, \dots, q \}$ choose indices $\mu, \nu \in \{ 1, \dots, p \}$
such that $\alpha_i^{} = \alpha_{\mu i}^{}$ and $\beta_i^{} = \beta_{\nu i}^{}$.
If $\mu = \nu$, then $\alpha_i^{} + \beta_i^{} = \alpha_{\mu i}^{} + \beta_{\mu i}^{} = \deg(f_{k_\mu}^{})_i^{} \le a_i^{}$.
If $\mu \neq \nu$, then
\[
\alpha_i^{} + \beta_i^{} = \alpha_{\mu i}^{} + \beta_{\nu i}^{} = \deg(\lcm(f_{k_\mu}^{}, f_{k_\nu}^{}))_i^{} \le a_i^{}.
\]
Since $i$ was chosen arbitrarily, we obtain the inequality
\[
\deg(\lcm(f_{k_1}^{}, \dots, f_{k_p}^{})) = (\alpha_1^{} + \beta_1^{}, \dots, \alpha_q^{} + \beta_q^{}) \le a.
\]
Thus, $(f_{k_1}^{}, \dots, f_{k_p}^{})$ must lie in $\Gamma_a^p(J)$.
This proves the reverse inclusion ``$\supset$''.
\end{proof}

\begin{proposition}
\label{vanishing}
Let $J \subset \KK[x_1^{}, y_1^{}, \dots, x_q^{}, y_q^{}]$ be a monomial ideal.
Consider $a \in \ZZ_+^q$ and let $\{ g_1^{}, \dots, g_n^{} \}$
be the set of minimal generators of $J$ whose degree is less or equal to $a$.
Assume that $\deg(\lcm(g_1^{}, g_l^{})) \le a$ for all indices $l \in \{ 2, \dots, n \}$.
We claim that $\der \hilb_{\SSS / J}^{}(a) = 0$.
\end{proposition}

\begin{proof}
Recall the sets $\Gamma_a^p(J)$ introduced above Lemma~\ref{delta}.
By hypothesis, $(g_1^{}, g_l^{})$ lies in $\Gamma_a^2(J)$ for all indices $l \in \{ 2, \dots, n \}$.
In view of Lemma~\ref{gammas},
for $2 \le p \le n$ we can write $\Gamma_a^p(J)$ as a disjoint union $\Gamma_a^p(J) = \Phi^p \sqcup \Psi^p$, where
\[
\Phi^p = \{ (g_1^{}, g_{l_2}^{}, \dots, g_{l_p}^{}) \mid 2 \le l_2^{} < \dots < l_p^{} \le n, \
(g_{l_2}^{}, \dots, g_{l_p}^{}) \in \Gamma_a^{p - 1}(J) \} \quad \text{and}
\]
\[
\Psi^p = \{ (g_{l_1}^{}, \dots, g_{l_p}^{}) \in \Gamma_a^p(J) \mid 2 \le l_1^{} < \dots < l_p^{} \le n \}.
\]
Notice that $\left| \Phi^p \right| = \left| \Psi^{p - 1} \right|$, where, by convention, $\Psi^1 = \{ g_2^{}, \dots, g_n^{} \}$.
Also notice that $\Psi^n = \emptyset$.
Applying Lemma~\ref{delta}, we calculate:
\begin{align*}
\der \hilb_{\SSS / J}^{}(a) & = 1 + \sum_{1 \le p \le n} (-1)^p \left| \Gamma_a^p(J) \right| \\
& = 1 - \left| \Gamma_a^1(J) \right| + \sum_{2 \le p \le n} (-1)^p (\left| \Phi^p \right| + \left| \Psi^p \right|) \\
& = 1 - n + \sum_{2 \le p \le n} (-1)^p \left| \Psi^{p - 1} \right| + \sum_{2 \le p \le n - 1} (-1)^p \left| \Psi^p \right| \\
& = 1 - n + \left| \Psi^1 \right| + \sum_{2 \le p \le n - 1} (-1)^{p + 1} \left| \Psi^p \right| + \sum_{2 \le p \le n - 1} (-1)^p \left| \Psi^p \right| \\
& = 1 - n + (n - 1) + \sum_{2 \le p \le n - 1} ((-1)^{p + 1} + (-1)^p) \left| \Psi^p \right| \\
& = 0. \qedhere
\end{align*}
\end{proof}


\section{Finite subschemes of a product of projective lines}
\label{end}

\noindent
In this section we assume that $\KK$ is infinite and that $\VV = (\PP^1)^q$, where $q \ge 2$.
We write $\SSS = \KK[x_1^{}, y_1^{}, \dots, x_q^{}, y_q^{}]$,
where $\deg(x_i^{}) = e_i^{}$ and $\deg(y_i^{}) = e_i^{}$.
Let $X \subset \VV$ be a zero-dimensional subscheme.
We recall from section~\ref{rectangular}
that $s_i^{}$ denotes the length of the projection $Z_i^{}$ of $X$ onto the $i$-th component of $\VV$.
Our first goal is to give a different proof for the fact
that the domain $[0, s_1^{} - 1] \times \dots \times [0, s_q^{} - 1]$ is relevant to $\hilb_X^{}$.
Our second goal is to give a procedure
for detecting quasi-rectangular domains that are relevant to $\hilb_X^{}$.
All results of this section are applications of Proposition~\ref{vanishing} and of Theorem~\ref{macaulay}.

\begin{lemma}
\label{lexicographic}
We adopt the above notations.
We consider the lexicographic monomial ordering on $\SSS$
such that $x_1^{} > y_1^{} > \dots > x_q^{} > y_q^{}$.
We assume that $y_1^{}$ does not vanish at any point of $\red(X)$.
We claim that $x_1^{s_1}$ is a minimal generator of $\lead(I(X))$ (see Notation~\ref{leader}).
We further claim that every other minimal generator of $\lead(I(X))$
has the form $x_1^{\alpha_1} x_2^{\alpha_2} y_2^{\beta_2} \cdots x_q^{\alpha_q} y_q^{\beta_q}$
with $0 \le \alpha_1^{} \le s_1^{} - 1$.
\end{lemma}

\begin{proof}
Write $J = \lead(I(X))$.
First we claim that $y_1^{}$ does not divide any minimal generator of $J$.
To prove this claim we argue by contradiction.
Assume that there existed a minimal generator $g$ of $J$ of the form
$g = x_1^{\alpha_1} y_1^{\beta_1} \cdots x_q^{\alpha_q} y_q^{\beta_q}$ with $\beta_1^{} > 0$.
Write $g = \lead(f)$ for some $\ZZ^q$-homogeneous polynomial $f \in I(X)$.
For any other monomial $u = x_1^{\gamma_1} y_1^{\delta_1} \cdots x_q^{\gamma_q} y_q^{\delta_q}$
occurring in $f$ we have the inequality $\alpha_1^{} \ge \gamma_1^{}$ because $g > u$
and we have the equation $\alpha_1^{} + \beta_1^{} = \gamma_1^{} + \delta_1^{}$ because $\deg(g) = \deg(u)$.
It follows that $\beta_1^{} \le \delta_1^{}$.
Since $u$ was chosen arbitrarily, it follows that $f$ is divisible by $y_1^{\beta_1}$.
Since $y_1^{}$ does not vanish at any point of $\red(X)$,
it follows that $f / y_1^{\beta_1}$ lies in $I(X)$.
Thus, $g / y_1^{\beta_1} = \lead(f / y_1^{\beta_1})$ belongs to $J$.
This contradicts the fact that $g$ is a minimal generator of $J$
and concludes the proof of the claim.

The ideal $J_1^{} = J \cap \KK[x_1^{}, y_1^{}]$ is a $\ZZ$-homogeneous ideal
of the $\ZZ$-graded ring $\SSS_1^{} = \KK[x_1^{}, y_1^{}]$.
According to Theorem~\ref{macaulay}, $\hilb_{\SSS / J}^{} = \hilb_X^{}$ hence
\begin{align*}
\hilb_{\SSS_1 / J_1}^{}(a_1^{}) & = \hilb_{\SSS / J}^{}(a_1^{}, 0, \dots, 0) = \hilb_X^{}(a_1^{}, 0, \dots, 0) = \hilb_{Z_1}^{}(a_1) \\
& = \begin{cases}
a_1^{} + 1 & \text{if $0 \le a_1^{} \le s_1^{} - 1$}, \\
s_1^{} & \text{if $a_1^{} \ge s_1^{}$}.
\end{cases}
\end{align*}
It follows that $J_1^{}$ is generated by a single monomial of degree $s_1^{}$,
which, according to the above claim, is not divisible by $y_1^{}$.
We deduce that $J_1^{}$ is generated by $x_1^{s_1}$.
This monomial must be a minimal generator of $J$.
Every other minimal generator of $J$ is not divisible by $y_1^{}$ or by $x_1^{s_1}$,
so it has the form given in the lemma.
\end{proof}

\noindent
As an application of our methods we obtain a third proof for a particular case of Theorem~\ref{main}.
We formulate this as a separate proposition.

\begin{proposition}
\label{lines_2}
Assume that $\KK$ is infinite.
Let $X \subset (\PP^1)^q$ be a zero-dimensional subscheme.
For each index $i \in \{ 1, \dots, q \}$,
let $s_i^{}$ be the length of the projection of $X$ onto the $i$-th copy of $\PP^1$.
We claim that $R(X) = [0, s_1^{} - 1] \times \dots \times [0, s_q^{} - 1]$
is the smallest rectangular relevant domain for $\hilb_X^{}$.
\end{proposition}

\begin{proof}
We must show that $\der \hilb_X^{}$ vanishes on the complement of $R(X)$,
i.e., we must show that  $\der \hilb_X^{}(a) = 0$ if $a_i^{} \ge s_i^{}$ for some index $i \in \{ 1, \dots, q \}$.
By symmetry, it is enough to consider only the case when $i = 1$.
Performing, if necessary, a linear change of coordinates on the first copy of $\PP^1$,
we may assume that $y_1^{}$ does not vanish at any point of $\red(X)$.
We choose a monomial ordering on $\SSS$ as in Lemma~\ref{lexicographic}
and we consider the ideal $J = \lead(I(X))$ (see Notation~\ref{leader}).
According to Theorem~\ref{macaulay}, $\hilb_X^{} = \hilb_{\SSS / J}^{}$,
so the theorem reduces to showing that $\der \hilb_{\SSS / J}^{}(a) = 0$ if $a_1^{} \ge s_1^{}$ and $a \ge 0$.
This will follow from Proposition~\ref{vanishing}.

We now verify the hypotheses of Proposition~\ref{vanishing}.
Consider the set $\{ g_1^{}, \dots, g_n^{} \}$ of minimal generators of $J$
whose degree is less or equal to $a$.
According to Lemma~\ref{lexicographic}, $x_1^{s_1}$ is a minimal generator of $J$.
By hypothesis $\deg(x_1^{s_1}) = (s_1^{}, 0, \dots, 0) \le a$, so we may take $g_1^{} = x_1^{s_1}$.
Applying again Lemma~\ref{lexicographic}, we see that for each index $l \in \{ 2, \dots, n \}$
we may write $g_l^{} = x_1^{\alpha_1} x_2^{\alpha_2} y_2^{\beta_2} \cdots x_q^{\alpha_q} y_q^{\beta_q}$
with $0 \le \alpha_1^{} \le s_1^{} - 1$.
Thus, $\lcm(g_1^{}, g_l^{}) = x_1^{s_1} x_2^{\alpha_2} y_2^{\beta_2} \cdots x_q^{\alpha_q} y_q^{\beta_q}$.
We have $\deg(\lcm(g_1^{}, g_l^{})) \le a$ because $s_1^{} \le a_1^{}$
and $\alpha_j^{} + \beta_j^{} \le a_j^{}$ for all indices $j \in \{ 2, \dots, q \}$,
by virtue of the fact that $\deg(g_l^{}) \le a$.
Thus, the hypotheses of Proposition~\ref{vanishing} are satisfied and we conclude that $\der \hilb_{\SSS / J}^{}(a) = 0$.
\end{proof}

\begin{remark}
\label{generators}
Let $I \subset \SSS$ be a $\ZZ^q$-homogeneous ideal.
Let $F$ and $G$ be two finite sets of generators of $I$ consisting of $\ZZ^q$-homogeneous non-zero polynomials.
We assume that $F$ is minimal, i.e.\ no proper subset of $F$ can generate $I$.
We assert that for every $f \in F$ there is $g \in G$ such that $\deg(f) = \deg(g)$.
Indeed, write $f = \sum_{g \in G} v_g^{} g$.
For each $g \in G$ write $g = \sum_{h \in F} w_{gh}^{} h$.
We may assume that all $v_g^{}$ and $w_{gh}^{}$ are $\ZZ^q$-homogeneous.
We claim that there is $g \in G$ such that $v_g^{} \neq 0$ and $w_{gf}^{} \neq 0$.
If this were not the case,
then $f$ would be a combination of elements in $F \setminus \{ f \}$,
which would contradict the minimality of $F$.
We have the relations $\deg(f) = \deg(v_g^{}) + \deg(g)$ and $\deg(g) = \deg(w_{gf}) + \deg(f)$,
hence $\deg(f) \ge \deg(g) \ge \deg(f)$.
\end{remark}

\begin{lemma}
\label{monomial}
Let $X \subset (\PP^1)^q$ be a zero-dimensional subscheme.
We assume that $I(X)$ is a monomial ideal.
We claim that the set of minimal generators of $I(X)$ is of the form $\{ u_1^{}, \dots, u_q^{}, g_1^{}, \dots, g_n^{} \}$,
where $n \ge 0$, $u_i^{} = x_i^{\alpha_i} y_i^{s_i - \alpha_i}$ for all $i = 1, \dots, q$,
and $\deg(g_l^{}) \le (s_1^{} - 1, \dots, s_q^{} - 1)$ for all $l = 1, \dots, n$.
\end{lemma}

\begin{proof}
Let $u_i^{}$ be the generator of $I(X) \cap \KK[x_i^{}, y_i^{}]$
and let $g_1^{}, \dots, g_n^{}$ be the minimal generators of $I(X)$
that do not lie in any $\KK[x_i^{}, y_i^{}]$.
We concentrate on proving the inequalities $\deg(g_l^{}) \le (s_1^{} - 1, \dots, s_q^{} - 1)$,
the rest of the lemma being obvious.
Since $\KK$ is infinite,
we can find $\kappa_1^{} \in \KK \setminus \{ 0 \}$
such that $z_1^{} = \kappa_1^{} x_1^{} + y_1^{}$ does not vanish at any point of $\red(X)$.
Regarding $\SSS$ as a polynomial ring in the variables $x_1^{}, z_1^{}, x_2^{}, y_2^{}, \dots, x_q^{}, y_q^{}$,
we consider the lexicographic ordering on $\SSS$
such that $x_1^{} > z_1^{} > x_2^{} > y_2^{} > \dots > x_q^{} > y_q^{}$.
Put $J = \lead(I(X))$ (see Notation~\ref{leader}).
According to Lemma~\ref{lexicographic}, $x_1^{s_1} = \lead(u_1^{})$ is a minimal generator of $J$
and every other minimal generator $v$ of $J$ satisfies the condition $\deg(v)_1^{} \le s_1^{} - 1$.
Let $G$ be a Gr\"obner basis of $I(X)$ containing $u_1^{}$ and consisting of $\ZZ^q$-homogeneous polynomials,
such that $\lead(G)$ is the set of minimal generators of $J$.
For every $g \in G \setminus \{ u_1^{} \}$ we have the relations $\deg(g)_1^{} = \deg(\lead(g))_1^{} \le s_1^{} - 1$.
We apply Remark~\ref{generators} to the set $F$ of minimal generators of $I(X)$ and to $G$.
For each $g_l^{}$ there is $g \in G$ such that $\deg(g_l^{}) = \deg(g)$.
Since $g_l^{} \notin \KK[x_1^{}, y_1^{}]$, it follows that $g \neq u_1^{}$,
hence $\deg(g_l^{})_1^{} = \deg(g)_1^{} \le s_1^{} - 1$.

In the same manner, for all indices $l \in \{ 1, \dots, n \}$ and $i \in \{1, \dots, q \}$,
by replacing the variable $y_i^{}$ with a suitable variable $z_i^{} = \kappa_i^{} x_i^{} + y_i^{}$,
chosen so as not to vanish at any point of $\red(X)$,
we can prove the inequality $\deg(g_l^{})_i^{} \le s_i^{} - 1$.
\end{proof}

\noindent
Let $X$ be as in Lemma~\ref{monomial}.
We recall, from Proposition~\ref{regularity_3}(v), the rectangular relevant domain
\[
R(X) = [0, s_1^{} - 1] \times \dots \times [0, s_q^{} - 1] = [0, \rem(X)] \subset \ZZ^q.
\]
When $n = 1$ we put $R_1^{} = [\deg(g_1^{}), \rem(X)]$.
When $n > 1$, for each index $l \in \{1, \dots, n \}$, we put
\[
R_l^{} =  \bigcap_{k \in \{ 1, \dots, n \} \setminus \{ l \}} [\deg(\lcm(g_k^{}, g_l^{})), \, \rem(X)].
\]

\begin{proposition}
\label{procedure}
Assume that $\KK$ is infinite.
Let $X \subset (\PP^1)^q$ be a zero-dimensional subscheme.
Assume that $I(X)$ is a monomial ideal
and that $X$ is not a complete intersection.
We claim that the quasi-rectangular domain
\[
D(X) = R(X) \setminus \bigcup_{1 \le l \le n} R_l^{} \subset \ZZ^q
\]
is a relevant domain for $\hilb_X^{}$.
If $I(X)$ has $q + 1$ minimal generators, then $D(X)$ is strictly contained in $R(X)$.
\end{proposition}

\begin{proof}
Under the notations of Lemma~\ref{monomial},
the hypothesis that $X$ not be a complete intersection is equivalent to saying that $n \ge 1$.
We must show that $\der \hilb_X^{}(a) = 0$ for all $a \in \ZZ^q \setminus D(X)$.
We already know from Theorem~\ref{main}
that $R(X)$ is a relevant domain for $\hilb_X{}$,
hence we may assume that $a \in R(X)$, i.e.\ $a \in R_l^{}$ for some index $l \in \{ 1, \dots, n \}$.
Relabeling $\{ g_1^{}, \dots, g_n^{} \}$, if necessary, we may take $l = 1$.
We shall apply Proposition~\ref{vanishing} with $J = I(X)$.
We now verify the hypotheses of Proposition~\ref{vanishing}.
By the construction of $R_1^{}$,
$a \ge \deg(\lcm(g_1^{}, g_k^{})) \ge \deg(g_k^{})$ for all indices $k \in \{ 2, \dots, n \}$ and $a \ge \deg(g_1^{})$.
Since $a$ belongs to $R(X)$, $a \ngeq \deg(u_i^{})$ for all indices $i \in \{ 1, \dots, q \}$.
From Lemma~\ref{monomial} we deduce
that $\{ g_1^{}, \dots, g_n^{} \}$ is, indeed, the set of minimal generators of $I(X)$
whose degree is less or equal to $a$.
The inequality from Proposition~\ref{vanishing} is satisfied by the definition of $R_1^{}$.
Thus, the hypotheses of Proposition~\ref{vanishing} are satisfied
and we conclude that $\der \hilb_X^{}(a) = 0$.

If $I(X)$ has $q + 1$ minimal generators, that is, if $n = 1$,
then, in view of Lemma~\ref{monomial}, $R_1^{} \neq \emptyset$,
forcing $D(X)$ to be strictly contained in $R(X)$.
\end{proof}

\noindent
If $n > 1$, then $R_l^{}$ may be empty for all indices $l$,
i.e.\ $D(X)$ may coincide with $R(X)$.
We finish this section with two examples in which $n > 1$, $X$ is non-ACM
and $D(X)$ is strictly contained in $R(X)$.

\begin{example}
\label{quadric_2}
Take $X \subset \PP^1(\CC) \times \PP^1(\CC)$ to be the union of two multiple points with ideals
$(x_1^{\alpha_1}, x_1^{} x_2^{}, x_2^{\alpha_2})$, respectively, $(y_1^{\beta_1}, y_1^{} y_2^{}, y_2^{\beta_2})$.
Here $\alpha_1^{}, \alpha_2^{}, \beta_1^{}, \beta_2^{} \ge 2$.
We have
\[
I(X) = (x_1^{\alpha_1} y_1^{\beta_1}, \, x_2^{\alpha_2} y_2^{\beta_2}, \, g_1^{}, \dots, g_7^{}),
\]
where
\begin{align*}
g_1^{} & = x_1^{\alpha_1} y_1^{} y_2^{},
& g_2^{} & = x_1^{\alpha_1} y_2^{\beta_2},
& g_3^{} & = x_1^{} y_1^{\beta_1} x_2^{},
& g_4^{} &= x_1^{} y_1^{} x_2^{} y_2^{}, \\
g_5^{} & = x_1^{} x_2^{} y_2^{\beta_2},
& g_6^{} & = y_1^{\beta_1} x_2^{\alpha_2},
& g_7^{} & = y_1^{} x_2^{\alpha_2} y_2^{}.
\end{align*}
We have the equations $s_1^{} = \alpha_1^{} + \beta_1^{}$ and $s_2^{} = \alpha_2^{} + \beta_2^{}$.
We have the relations
\begin{align*}
\lcm(g_1^{}, g_3^{}) & = x_1^{\alpha_1} y_1^{\beta_1} x_2^{} y_2^{},
& \deg(\lcm(g_1^{}, g_3^{})) & = (\alpha_1^{} + \beta_1^{}, 2) \nleq (s_1^{} - 1, s_2^{} - 1), \\
\lcm(g_2^{}, g_6^{}) & = x_1^{\alpha_1} y_1^{\beta_1} x_2^{\alpha_2} y_2^{\beta_2},
& \deg(\lcm(g_2^{}, g_6^{})) & = (\alpha_1^{} + \beta_1^{}, \alpha_2^{} + \beta_2^{}) \nleq (s_1^{} - 1, s_2^{} - 1), \\
\lcm(g_5^{}, g_7^{}) & = x_1^{} y_1^{} x_2^{\alpha_2} y_2^{\beta_2},
& \deg(\lcm(g_5^{}, g_7^{})) & = (2, \alpha_2^{} + \beta_2^{}) \nleq (s_1^{} - 1, s_2^{} - 1).
\end{align*}
We deduce that $R_1^{}$, $R_3^{}$, $R_2^{}$, $R_6^{}$, $R_5^{}$ and $R_7^{}$ are empty.
We have the equations
\begin{align*}
\lcm(g_4^{}, g_1^{}) & = x_1^{\alpha_1} y_1^{} x_2^{} y_2^{},
& \lcm(g_4^{}, g_2^{}) & = x_1^{\alpha_1} y_1^{} x_2^{} y_2^{\beta_2}, \\
\lcm(g_4^{}, g_3^{}) & = x_1^{} y_1^{\beta_1} x_2^{} y_2^{},
& \lcm(g_4^{}, g_5^{}) & = x_1^{} y_1^{} x_2^{} y_2^{\beta_2}, \\
\lcm(g_4^{}, g_6^{}) & = x_1^{} y_1^{\beta_1} x_2^{\alpha_2} y_2^{},
& \lcm(g_4^{}, g_7^{}) & = x_1^{} y_1^{} x_2^{\alpha_2} y_2^{}.
\end{align*}
Thus,
\begin{align*}
\max_{k = 1, 2, 3, 5, 6, 7} \deg(\lcm(g_4^{}, g_k^{}))_1^{}
& = 1 + \max \{ \alpha_1^{}, \beta_1^{} \} \le s_1^{} - 1, \\
\max_{k = 1, 2, 3, 5, 6, 7} \deg(\lcm(g_4^{}, g_k^{}))_2^{}
& = 1 + \max \{ \alpha_2^{}, \beta_2^{} \} \le s_2^{} - 1.
\end{align*}
We deduce that $R_4^{} \neq \emptyset$.
In fact,
\[
R_4^{} = [1 + \max \{ \alpha_1^{}, \beta_1^{} \}, \, \alpha_1^{} + \beta_1^{} - 1] \times
[1 + \max \{ \alpha_2^{}, \beta_2^{} \}, \, \alpha_2^{} + \beta_2^{} - 1].
\]
We conclude that
\[
D(X) = [0, \, \alpha_1^{} + \beta_1^{} - 1] \times [0, \, \alpha_2^{} + \beta_2^{} - 1] \setminus R_4^{}.
\]
If $X$ were ACM, then, in view of Theorem~\ref{giuffrida}(iii),
the degrees of $g_1^{}, \dots, g_7^{}$ would be incomparable.
However, $\deg(g_4^{}) = (2, 2) \le \deg(g_2^{}) = (\alpha_1^{}, \beta_2^{})$.
This shows that $X$ is not ACM.
\end{example}

\begin{example}
\label{quadric_3}
Take $X \subset \PP^1(\CC) \times \PP^1(\CC)$ to be the union of three multiple points with ideals
$(x_1^{\alpha_1}, x_2^{\alpha_2})$, $(x_1^\alpha, y_2^\beta)$, respectively, $(y_1^{\beta_1}, y_2^{\beta_2})$.
We assume that $\alpha_1^{} < \alpha$ and $\beta < \beta_2^{}$.
We have
\[
I(X) = (x_1^\alpha y_1^{\beta_1}, \, x_2^{\alpha_2} y_2^{\beta_2}, \, g_1^{}, \, g_2^{}, \, g_3^{}),
\]
where $g_1^{} = x_1^{\alpha_1} y_1^{\beta_1} y_2^\beta$, $g_2^{} = x_1^{\alpha_1} y_2^{\beta_2}$ and
$g_3^{} = y_1^{\beta_1} x_2^{\alpha_2} y_2^\beta$.
We have the equations $s_1^{} = \alpha + \beta_1^{}$ and $s_2^{} = \alpha_2^{} + \beta_2^{}$.
We have the relations
\begin{align*}
\lcm(g_1^{}, g_2^{}) & = x_1^{\alpha_1} y_1^{\beta_1} y_2^{\beta_2},
& \deg(\lcm(g_1^{}, g_2^{})) & = (\alpha_1^{} + \beta_1^{}, \beta_2^{}) \le (s_1^{} - 1, s_2^{} - 1), \\
\lcm(g_1^{}, g_3^{}) & = x_1^{\alpha_1} y_1^{\beta_1} x_2^{\alpha_2} y_2^\beta,
& \deg(\lcm(g_1^{}, g_3^{})) & = (\alpha_1^{} + \beta_1^{}, \alpha_2^{} + \beta) \le (s_1^{} - 1, s_2^{} - 1), \\
\lcm(g_2^{}, g_3^{}) & = x_1^{\alpha_1} y_1^{\beta_1} x_2^{\alpha_2} y_2^{\beta_2},
& \deg(\lcm(g_2^{}, g_3^{})) & = (\alpha_1^{} + \beta_1^{}, \alpha_2^{} + \beta_2^{}) \nleq (s_1^{} - 1, s_2^{} - 1).
\end{align*}
We deduce that $R_2^{} = \emptyset$, $R_3^{} = \emptyset$ and that
\[
R_1^{} = [\alpha_1^{} + \beta_1^{}, \, \alpha + \beta_1^{} - 1] \times [\max \{ \beta_2^{}, \alpha_2^{} + \beta \}, \, \alpha_2^{} + \beta_2^{} - 1].
\]
We conclude that
\[
D(X) = [0, \, \alpha + \beta_1^{} - 1] \times [0, \, \alpha_2^{} + \beta_2^{} - 1] \setminus R_1^{}.
\]
Assume that $X$ were ACM.
Then, according to Theorem~\ref{giuffrida}(iii), there would exist homogeneous polynomials
$u_1^{}$, $u_2^{}$, $u_3^{}$, $u_4^{} \in \CC[x_1^{}, y_1^{}]$ and
$v_1^{}$, $v_2^{}$, $v_3^{}$, $v_4^{} \in \CC[x_2^{}, y_2^{}]$
such that
\[
I(X) = (u_1^{} u_2^{} u_3^{} u_4^{}, \ v_1^{} u_2^{} u_3^{} u_4^{}, \ v_1^{} v_2^{} u_3^{} u_4^{}, \ v_1^{} v_2^{} v_3^{} u_4^{}, \
v_1^{} v_2^{} v_3^{} v_4^{}).
\]
We choose $\kappa$ and $\lambda$ in $\CC^*$ such that $v_i^{}(\kappa, \lambda) \neq 0$ for $i = 1, 2, 3, 4$.
We reduce the above equality of ideals modulo $(x_2^{} - \kappa, \ y_2^{} - \lambda)$.
We obtain the equality of ideals
\[
(1) = (u_1^{} u_2^{} u_3^{} u_4^{}, \ u_2^{} u_3^{} u_4^{}, \ u_3^{} u_4^{}, \ u_4^{})
\]
in $\CC[x_1^{}, y_1^{}]$.
This is absurd.
We conclude that $X$ is not ACM.
\end{example}

\end{document}